\documentclass[preprint,12pt]{article}
\usepackage{amssymb}
\usepackage{mathrsfs}
\usepackage{amsfonts}
\usepackage{graphicx}
\usepackage{indentfirst}
\usepackage{colortbl,dcolumn}
\usepackage{color}
\usepackage{amsmath}
\usepackage{psfrag}
\usepackage{booktabs}
\usepackage{cite}
\usepackage{amsthm}
\usepackage{float}
\usepackage{verbatim}
\usepackage{appendix}
\numberwithin{equation}
{section}
\usepackage[top=1in, bottom=1in, left=0.8in, right=0.8in]{geometry}

\newcommand{\R}{\mathbb{R}}
\newcommand{\N}{\mathbb{N}}
\newcommand{\E}{\mathbb{E}}
\newcommand{\I}{\mathbb{I}}
\newcommand{\dd}{\text{d}}
\newtheorem{thm}{Theorem}[section]
\newtheorem{defn}[thm]{Definition}
\newtheorem{lem}[thm]{Lemma}

\newtheorem{cor}[thm]{Corollary}

\newtheorem{assumption}[thm]{Assumption}
\newcommand{\yuebox}[1]{\fbox{$\triangleright$\textcolor{blue}{\textbf{Yue}:} #1}}

\begin{document}

\title{{\color{black}A projected} Euler Method for Random Periodic Solutions of Semi-linear SDEs {\color{black}with non-globally Lipschitz coefficients}
\footnotemark[2] \footnotetext[2]{This work was supported by Natural Science Foundation of China (12071488, 12371417, 11971488). YW would like to acknowledge the support of the Royal Society through the International Exchanges scheme IES\textbackslash R3\textbackslash 233115. 
The authors want to thank Yinghao Huang for his suggestions and comments improving the manuscript.
                \\
                E-mail addresses:
 y.j.guo@foxmail.com, x.j.wang7@csu.edu.cn, 
yue.wu@strath.ac.uk. 
                }}
\author{Yujia Guo$^{a}$, Xiaojie Wang$^{a}$, and Yue Wu$^{b}$
\\
\footnotesize $^a$ School of Mathematics and Statistics, HNP-LAMA, Central South University, Changsha, \\
\footnotesize Hunan, P. R. China \\
\footnotesize $^b$ Department of Mathematics and Statistics, University of Strathclyde, Glasgow G1 1XH, UK} 
\date{\today}

\maketitle
\begin{abstract}
{\color{black}
The present work introduces and investigates
an explicit time discretization scheme, called the projected Euler method,
to numerically approximate  random periodic solutions of semi-linear SDEs  under
non-globally Lipschitz conditions.
The existence of the random periodic solution is demonstrated as the limit of the pull-back of the discretized SDE.
Without relying on a priori high-order moment bounds of the numerical approximations, the mean square convergence rate of the approximation scheme is proved to be
order $0.5$ for SDEs with multiplicative noise and order $1$ for SDEs with additive noise.
Numerical examples are also provided to validate our theoretical findings.}\\


\end{abstract}

\section{Introduction}
Periodic occurrences abound throughout nature. Since the pioneering works of Poincar\'e \cite{poincare1893methodes}, periodicity has consistently remained a focal point in the examination of dynamical systems. It has garnered significant interest in various fields including thermodynamics \cite{prodan2003kohn}, porous media \cite{adler2010resurgence}, quantum time crystals \cite{nakatsugawa2019time}, Thomas-Fermi plasma \cite{ur2020periodic}, and numerous other domains. Nonetheless, many real-world issues are prone to random fluctuations induced by uncertainty and unknown variables. Hence, the exploration of random periodicity emerges as a fundamentally crucial area of study.

As the stochastic counterpart to periodic solutions, the definition of random periodic solutions for a $C^{1}$-cocycle was initially proposed by Zhao and Feng \cite{zhao2009random}, while Feng, Zhao, and Zhou \cite{feng2011pathwise} subsequently expanded upon this concept for semiflows. Their work has catalyzed further advancements in the exploration of various issues within autonomous and non-autonomous stochastic differential equations. This includes investigations into the existence of random solutions generated by non-autonomous SPDEs with additive noise \cite{feng2012random}, the anticipation of random solutions of SDEs with multiplicative linear noise \cite{feng2016anticipating}, periodic measures and ergodicity \cite{feng2020random}, among others
\cite{uda2024random,feng2021ergodic}.

Given 
$W:\R \times \Omega \rightarrow
\R^m$
a standard two-sided Wiener process on the 
complete
probability space $(\Omega,\mathcal{F},\mathbb{P})$, where
the filtration is defined as
$\mathcal{F}^{t}_{s}:=
\sigma\{W_{u}-W_{v}:s\leq v\leq u\leq t\}$ and $\mathcal{F}^{t}
=\mathcal{F}^{t}_{-\infty}
=\bigvee_{s\leq t}\mathcal{F}^{t}_{s}$.
We consider the  following semi-linear SDEs with multiplicative noise:

\begin{equation}
	\label{eq_PEM:Problem_SDE}
	\left\{
	\begin{aligned}
	   \dd X_{t}^{t_{0}} 
         & = 
         \big(
         A X_{t}^{t_{0}}
         +
         f(t, \,X_{t}^{t_{0}})
         \big)
         \dd t
         + 
         g(t,
         \,X_{t}^{t_{0}})
         \, \dd W_t,
   \quad t \in (t_{0},T],\\
		 X_{t_{0}}^{t_{0}} & = \xi,
	\end{aligned}\right.
\end{equation}
where 
$A \in \R^{d \times d}$ is a negative-definite matrix,
$f:\mathbb{R} \times \mathbb{R}^{d}\rightarrow \mathbb{R}^{d}$,  $g:\mathbb{R} \times \mathbb{R}^{d} \rightarrow  \mathbb{R}^{d \times m}$ are continuous functions.
We use $X_{t_1}^{t_0}$ to emphasise a process $X$ evaluated at $t_1$ which starts from $t_0$.
The random initial value $\xi$ is assumed to be $\mathcal{F}^{t_{0}}$-measurable. Note that by the variation of constant formula, the solution of \eqref{eq_PEM:Problem_SDE} can be written as
\begin{equation}
    \label{eq:variation of comstant formula}
    X^{t_{0}}_{t}(\xi)=
    e^{A (t-t_{0})}\xi
    +\int_{t_{0}}^{t}
    e^{A (t-s)} f(s,X^{t_{0}}_{s})
    \,\dd s
    +\int_{t_{0}}^{t}
    e^{A(t-s)} 
    g(s,X^{t_{0}}_{s})
    \,\dd W_{s}.
\end{equation}

The relevant research on the numeric of  SDEs has made rapid progress recently
\cite{hutzenthaler2020perturbation,hutzenthaler2015numerical,milstein2004stochastic,kloeden1992stochastic,li2019explicit}.
In general, the explicit computation of random periodic solutions is often unattainable, necessitating the utilization of numerical approximations, which play a pivotal role in this domain. The initial study by Feng et al. \cite{feng2017numerical} employed classical numerical methods, such as the Euler-Maruyama method and a modified Milstein method, to approximate random periodic solutions for a dissipative system with global Lipschitz conditions. Wei and Chen \cite{wei2020numerical} subsequently extended the applicability of the Euler-Maruyama method to the stochastic theta method, demonstrating convergence to the exact solution at an order of $1/4$.
{\color{black}Moradi et al.\cite{moradi2024random} further explored this topic by simulating random periodic solutions using $\theta$-Maruyama and $\theta$-Milstein methods with weaker conditions on the drift term. }


Wu \cite{wu2023backward} delved into the study of the existence and uniqueness of random periodic solutions for an additive SDE with a one-sided Lipschitz condition and provided an analysis indicating an order-half convergence of its numerical approximation using the backward Euler method. Later, Guo, Wang, and Wu \cite{guo2023order} {\color{black}lifted} the convergence order from half to one under a relaxed condition compared to \cite{wu2023backward}. 
Recently, Chen et al. \cite{chen2024stochastic}
turned to stochastic theta
methods and showed that the mean square convergence order is 0.5 for SDEs with multiplicative noise and 1 for SDEs with additive noise under non-globally Lipschitz conditions.

Different from works mentioned above, in this article we consider explicit time-stepping schemes for the numerical approximation to random periodic solution of semi-linear SDEs under non-globally Lipschitz conditions.
The conditions are weaker compared to literature \cite{wu2023backward,guo2023order}.
This applies the projected technique, previously used in \cite{beyn2016stochastic,beyn2017stochastic} for SDEs in finite time interval, to derive convergence results for random periodic solutions in infinite time intervals.
The projected Euler method involves the standard Euler method combined with a projection onto a ball that expands in radius  with a negative power of the step size.
This approach helps prevent the nonlinear drift and diffusion from causing  excessively large values, even in infinite time horizon.

The main focus of this paper is to analyze the strong convergence rate of the projected Euler method applied to the random periodic solution of semi-linear SDEs under non-global conditions. 
Without relying on a priori high-order moment bounds of the numerical approximations, we determine that the mean square convergence order is $0.5$ for SDEs with multiplicative noise and $1$ for SDEs with additive noise.

The paper is structured as follows: Section \ref{sec:Random Periodic Solutions of SDEs} outlines the standard notation and assumptions utilized in our proofs, and establishes the existence and uniqueness of the random periodic solution. In Section \ref{sec:PEM}, we detail the well-posedness and the existence of a unique random periodic solution using the projected Euler method. Section \ref{sec:strong_rate_of_PEM} is dedicated to the error analysis concerning random periodic solutions derived from the projected Euler method. Finally, Section \ref{sec:numerical results} presents several numerical experiments aimed at illustrating the theoretical findings.

\section{Random Periodic Solutions of SDEs}
\label{sec:Random Periodic Solutions of SDEs}
Recalling the definition of the random periodic solution for stochastic semi-flows  given in \cite{zhao2009random}.
Let $X$ be a separable Banach space.
Denote by
$ (\Omega,\mathcal{F},\mathbb{P},(\theta _{s})_{s \in \mathbb{R}}
)$ 
a metric dynamical system and 
$\theta _{s}:\Omega \rightarrow \Omega$ 
is assumed to be a measurably invertible 
for all $s \in \mathbb{R}$.
Denote $\Delta:= \{(t,s)  \in \mathbb{R}^{2},s \leq t \}$.  
Consider a stochastic  periodic semi-flow 
$u:\Delta \times \Omega \times X \rightarrow X $ of period $\tau$, 
which satisfies the following standard condition
\begin{equation}
    u(t,r,\omega) =
    u(t,s,\omega) 
    \circ
    u(s,r,\omega),
\end{equation}
and the periodic property
\begin{equation}
u(t+\tau,s+\tau,\omega)=
u(t,s,\theta_{\tau}\omega),
\end{equation}
for all $r \leq s \leq t,r,s,\in \mathbb{R}$, 
for a.e. $\omega \in \Omega$.

\begin{defn}\label{def:rps}
	A random periodic solution of period $\tau
       >0 $ of a semi-flow
	 $u:\Delta \times \Omega \times X \rightarrow X 
       $
	  is an $\mathcal{F}$-measurable map 
	  $Y: \mathbb{R} \times \Omega \rightarrow X $ 
	  such that
	\begin{equation}
		u(t+\tau,t,Y(t,\omega),\omega)=
		Y(t+\tau,\omega)=
		Y(t,\theta_{\tau}\omega).
	\end{equation}
	for any $(t,s) \in \Delta$, $\omega \in \Omega$.
\end{defn}

Throughout this paper the following notation
is frequently used.
For simplicity, we denote $[d]:=
\{1,...,d\}$ and the letter
$ C $ is used to denote a generic positive
constant independent of time step size and
may vary for each appearance. 
Let $|\cdot|$, $ \| \cdot \| $ and 
$ \left\langle \cdot , \cdot \right\rangle $
be the absolute value of a scalar, the Euclidean norm and the inner product of vectors
in $ \R^d $, respectively.
By $A^{T}$ we denote the transpose of vector or matrix.
Given a matrix $A$, we use $\| A \|:=\sqrt{trace(A^{T}A)}$ to denote the trace norm of $A$.
On a probability space 
$ ( \Omega, \mathcal{ F }, \mathbb{P} ) $, 
we use $ \E $ to denote the mean expectation and
$ L^p(\Omega; \textcolor{black}{\R^{d }}) $, 
$ d \in \mathbb{N} $,
to denote the family of 
$ \textcolor{black}{\R^{d}} $-valued variables with
the norm defined by 
$ \| \xi \|_{L^p(\Omega;\textcolor{black}{\R^{d}})} 
=(\E [\|\xi \|^p])^{\frac{1}{p}} < \infty $.

 
We present the following assumptions required to establish our main results.

\begin{assumption}
\label{ass_PEM}
Suppose the following conditions are satisfied.

(\romannumeral1)
$A$ is self-adjoint and negative definite  and
there exists a non-decreasing sequence $(\lambda_{i})_{i \in [d] } \subset \mathbb{R}$ 
of positive real numbers 
and an orthonormal basis $(e_{i})_{i \in [d]}$, 
such that 
$A e_{i} =-\lambda_{i} e_{i}$,  $i \in [d]$.
Moreover, one also obtains
\begin{equation}
\label{eq_PEM:lambda}
    \langle
    x,Ax
    \rangle
    \leq
    - \lambda_{1}
    \| x \|^{2},
    \qquad 
    \forall x \in \R^{d}.
\end{equation}

(\romannumeral2)
The drift coefficient functions $f$ and diffusion coefficient functions $g$ are continuous and periodic in time with period $\tau >0$, i.e.,
\begin{equation}
    f(t+\tau,x)=f(t,x),
    \quad
   g(t+\tau,x)=g(t,x),
    \qquad
    \forall
    x \in \R^{d},
    t \in \R.
\end{equation}

(\romannumeral3)
For some $p_{1} \in (1,\infty)$, there exists a constant $\alpha_{1 }< \lambda_1$ such that 
for and $x,y \in \R^d$
and
$t \in [0,\tau)$
\begin{equation}
\label{eq:coupled__momotoncity_condition}
    \langle 
    x-y , 
    f(t,x)-f(t,y)
    \rangle
    +
    \dfrac{2p_{1}-1}{2}
    \|
    g(t,x)-g(t,y)
    \|^{2}
    \leq 
    \alpha_{1} 
     \| x-y \|^{2}.
\end{equation}

(\romannumeral4 )
There exists some positive  constant 
{$\color{black}\gamma \in \Big[1,\frac{p_1+1}{2}\Big)$},
for $C_1, C >0$
such that
\begin{align}
\label{eq_PEM:the_esti_f(t,x)-f(t,y)}
     \| 
      f(t,x) -f(t,y)
     \|
     \vee
     \| 
      g(t,x) -g(t,y)
     \|
     &\leq
     {\color{black}C_1}
       (1+
       \|x\|^{\gamma-1}+\|y \|^{\gamma-1})
       \|x-y\|,
       \qquad
       \forall x,y \in \R^{d},\\
\label{eq_PEM:the_esti_f(t,x)-f(s,y)}
      \|f(t,x)-f(s,x)\|
      \vee
      \|g(t,x)-g(s,x)\|
      &\leq
       C(1+\|x\|^{\gamma})
       |t-s|,
       \qquad
       \forall x \in \R^{d},
       s,t \in [0,\tau).
\end{align}
(\romannumeral5 )
For {\color{black} any $p \geq 1$}, 
there exists a constant $C^{*}>0$ depending on $p$ such that 
$\E \big[\|\xi\|^{2p} \big]\leq C^{*}$.
\end{assumption}

The spatial regularity in \eqref{eq_PEM:the_esti_f(t,x)-f(t,y)} of Assumption \ref{ass_PEM}
 immediately implies,
 {\color{black}there exists an $C_2>0$},
\begin{align}
\label{eq_PEM:the_esti_f(t,x)}
    \|f(t,x)\|
    &\leq
    {\color{black}C_2}
    (1+\|x\|^{\gamma}),
    \qquad
    \forall x \in \R^{d},\\
\label{eq:the_esti_g(t,x)}
    \|
    g(t,x)
    \|
   & \leq
    C
    (1+\|x \|^{\gamma}),
    \qquad
    \forall x \in \R^{d}.
\end{align}

It can be verified that Assumption \ref{ass_PEM} leads to the following estimates.


\begin{lem}
\label{lem::coupled_monononcity_condition}
Let Assumption \ref{ass_PEM} be fulfilled,
for any $p_2 \in [1,p_1)$, there exists a small 
 positive constant $\epsilon$
such that 
\begin{equation}
\label{eq_PEM:coercivity_condition}
     \langle 
    x, 
    f(t,x)
    \rangle
    +
    \dfrac{2p_{2}-1}{2}
    \|
    g(t,x)
    \|^{2}
     \leq 
     \alpha_{2} 
     \|x\|^{2}
     +c_0,
     \qquad
     \forall
     x \in \R^{d},
     t \in \R,
\end{equation}
where $\alpha_2 
=\alpha_1+\epsilon<\lambda_1$,
$c_0=
\tfrac{\|f(t,0)\|^2}{2\epsilon}
+\tfrac{(2p_1-1)^2}{4(p_1-p_2)}\|g(t,0)\|^2
+\tfrac{2p_1-1}{2}\|g(t,0)\|^2$.
\end{lem}
The proof of 
Lemma \ref{lem::coupled_monononcity_condition} can be found in Appendix A.

The following assumption ensures the existence and uniqueness of a random periodic solution of \eqref{eq_PEM:Problem_SDE}
under non-globally Lipschitz conditions.
{\color{black}
\begin{assumption}
\label{thm:unique_random_periodic_solution}
Assume that
there exists a unique random periodic solution $X_{t}^{*}(\cdot) \in L^{2}(\Omega)$ with the form \begin{equation}
    \label{eq:pull-back}
    X^{*}_{t}
    =
    \int_{-\infty}^{t}
    e^{A(t-s)}
    f(s,X^{*}_{s})
    \,\dd s
    +\int_{-\infty}^{t}
    e^{A(t-s)}
    g(s,X^{*}_{s})
    \,\dd W_{s},
\end{equation}
such that $X^{*}$ is a limit of the pull-back 
$X^{-k\tau}_{t}(\xi)$ of \eqref{eq_PEM:Problem_SDE} when $k \rightarrow \infty$, i.e., 
	\begin{equation}
		\lim_{k \rightarrow \infty }
        \E\Big[
        \big\| X_{t}^{-k\tau}(\xi) -X_{t}^{*} \big\|^{2}
        \Big]
        =0.
	\end{equation}
\end{assumption}
}

Before moving on, we introduce a useful lemma for later use,
which has given in 
\cite{MR2016992}.
\begin{lem}
   \label{lem:Gronwall}
    Let $u,v,m$ be real-valued continuous functions defined on $[a,b]$,
    $ m(t) \geq 0$
    for $t \in [a,b]$.
    If $u$ satisfies the following inequality    
    \begin{equation}
        u(t)
        \leq
        v(t)
        +
        \int_{s}^{t}
        m(s)u(s)
        \, \mathrm{d} s,
    \end{equation}
    then
    \begin{equation}\label{eq:Gronwall 1}
        u(t)
        \leq
        v(t)
        +
        \int_{a}^{t}
        m(s)v(s)
        \exp
        \Big(
        \int_{s}^{t}
        m(r) \,\mathrm{d} r
        \Big)
        \, \mathrm{d} s.
    \end{equation}
    If in addition, the function $v$ is constant,
    then from
    \begin{equation}
        u(t)
        \leq
        v
        +
        \int_{s}^{t}
        m(s)u(s)
        \, \mathrm{d} s,
    \end{equation}
    it follows that
    \begin{equation}
    \label{eq:Gronwall 2}
        u(t)
        \leq
        v
        \exp
        \Big(
        \int_{a}^{t}
        m(r) \, \mathrm{d} r
        \Big).
    \end{equation}
\end{lem}

We first analyze the boundedness of the uniform moment of its solution under above assumptions.

\begin{lem}
\label{lem_PEM:the_pth_of_exact_solution}
    Let Assumption \ref{ass_PEM} hold, 
    consider the solution $X^{-k\tau}_{t}$ of SDE \eqref{eq_PEM:Problem_SDE}. If the initial value $X^{-k\tau}_{-k\tau}=\xi$,
    then for any {\color{black}$p
    \in [1,p_1)$}, there exists a positive constant $C$ 
    {\color{black}  depends on p}
    such that 
\begin{equation}\label{eq:the_p_th_(X(t))}
     \sup_{t\geq -k\tau}\E
     \Big[ 
     {\big\| X^{-k\tau}_{t} \big\|}^{2p} 
     \Big]
     \leq
      C
      <
      \infty.
\end{equation}
\end{lem}
\begin{proof}
    [Proof of Lemma \ref{lem_PEM:the_pth_of_exact_solution}.]
    {\color{black}
    Applying the It\^o formula to the following quantity for some constant $\epsilon_1>0$,
    \begin{equation}
    \begin{split}
        e^{2\lambda_{1}p
        (t+k\tau)}
        \big(
        \epsilon_1
        + {\big\| X^{-k\tau}_{t}\big\|}^{2}
        \big)^{p}
	    = &
		\big(
             \epsilon_1
             + {\|\xi\|}^{2}
            \big)^{p}
	    +
            2\lambda_1 p
           \int_{-k\tau}^{t}
            e^{2\lambda_{1}p(s+k\tau)}
            \big(
            \epsilon_1
            + {\big\| X^{-k\tau}_{s}\big\|}^{2}
            \big)^{p}
            \,\dd s \\
            & +
            2p
            \int_{-k\tau}^{t}
            e^{2\lambda_{1}p(s+k\tau)}
            \big
            (\epsilon_1
            +{\big\| X^{-k\tau}_{s}\big\|}^{2}
            \big)^{p-1}
            \big\langle 
             X^{-k\tau}_{s},
             A X^{-k\tau}_{s}
            \big\rangle
            \, \dd s \\
		&
		+2p \int_{-k\tau}^{t}
            e^{2\lambda_{1}p(s+k\tau)}
            \big( 
            \epsilon_1
            + {\big\| X^{-k\tau}_{s}\big\|}^{2} \big)^{p-1}
            \big\langle 
              X^{-k\tau}_{s},
              f(s,X^{-k\tau}_{s})
            \big\rangle
            \, \dd s 	\\	
		&
		+2p
            \int_{-k\tau}^{t}
            e^{2\lambda_{1}p(s+k\tau)}
            \big( 
            \epsilon_1
            +{\big\| X^{-k\tau}_{s}\big\|}^{2} \big)^{p-1}
		\big\langle 
		X^{-k\tau}_{s},
		g(s,X^{-k\tau}_{s})  
		\, \dd W_{s}
		\big\rangle \\
		&
		+p
            \int_{-k\tau}^{t}
            e^{2\lambda_{1}p(s+k\tau)}
            \big(
            \epsilon_1
            +{\big\|X^{-k\tau}_{s}\big\|}^{2} \big)^{p-1}
            \big\| 
            g(s,X^{-k\tau}_{s}) 
            \big\|^{2}
            \, \dd s\\
            & +
            2p(p-1)
            \int_{-k\tau}^{t}
            e^{2\lambda_{1}p(s+k\tau)}
            \big(
            \epsilon_1
            +{\big\|X^{-k\tau}_{s}\big\|}^{2} \big)^{p-2}
            \big\| 
            (X^{-k\tau}_s)^T
            g(s,X^{-k\tau}_{s}) 
            \big\|^{2}
            \, \dd s.
    \end{split}
    \end{equation}
    }
    {\color{black}Combining the last two terms on the right-hand-side gives}
\begin{equation}\label{eqn:ito1}
\begin{split}
        e^{2\lambda_{1}p
        (t+k\tau)}
        \big(
        \epsilon_1
        + {\big\| X^{-k\tau}_{t}\big\|}^{2}
        \big)^{p}
	    \leq &
		\big(
             \epsilon_1
             + {\|\xi\|}^{2}
            \big)^{p}
	    +
            2\lambda_1 p
           \int_{-k\tau}^{t}
            e^{2\lambda_{1}p(s+k\tau)}
            \big(
            \epsilon_1
            + {\big\| X^{-k\tau}_{s}\big\|}^{2}
            \big)^{p}
            \,\dd s \\
            & +
            2p
            \int_{-k\tau}^{t}
            e^{2\lambda_{1}p(s+k\tau)}
            \big
            (\epsilon_1
            +{\big\| X^{-k\tau}_{s}\big\|}^{2}
            \big)^{p-1}
            \big\langle 
             X^{-k\tau}_{s},
             A X^{-k\tau}_{s}
            \big\rangle
            \, \dd s \\
		&
		+2p \int_{-k\tau}^{t}
            e^{2\lambda_{1}p(s+k\tau)}
            \big( 
            \epsilon_1
            + {\big\| X^{-k\tau}_{s}\big\|}^{2} \big)^{p-1}
            \big\langle 
              X^{-k\tau}_{s},
              f(s,X^{-k\tau}_{s})
            \big\rangle
            \, \dd s 	\\	
		&
		+2p
            \int_{-k\tau}^{t}
            e^{2\lambda_{1}p(s+k\tau)}
            \big( 
            \epsilon_1
            +{\big\| X^{-k\tau}_{s}\big\|}^{2} \big)^{p-1}
		\big\langle 
		X^{-k\tau}_{s},
		g(s,X^{-k\tau}_{s})  
		\, \dd W_{s}
		\big\rangle \\
		&
		+p(2p-1)
            \int_{-k\tau}^{t}
            e^{2\lambda_{1}p(s+k\tau)}
            \big(
            \epsilon_1
            +{\big\|X^{-k\tau}_{s}\big\|}^{2} \big)^{p-1}
            \big\| 
            g(s,X^{-k\tau}_{s}) 
            \big\|^{2}
            \, \dd s.
\end{split}
\end{equation}
For every integers $n \geq 1$, define the stopping time
\begin{equation}
    \tau_{n}:= 
    \inf 
    \Big\{
       s \in [-k\tau,\infty):
      \big\|X_{s}^{-k\tau}\big\|
      \geq n
    \Big\}.
\end{equation}
Taking expectations on both sides of \eqref{eqn:ito1}, using \eqref{eq_PEM:lambda} of Assumption \ref{ass_PEM} 
and \eqref{eq_PEM:coercivity_condition} 
and letting $\epsilon_1 \rightarrow 0^+$ yield
\begin{equation}
\begin{split}
        & 
        \E
        \Big[
        e^{2\lambda_{1}p
        (t \wedge \tau_{n}+k\tau)}
        {\big\| X^{-k\tau}_{t \wedge \tau_{n}}\big\|}^{2p}
        \Big] \\
        & \leq 
         \E
         [ 
         \| \xi\| ^{2p}
         ]
         +
         \underbrace
         {
         2p
         \E 
         \Bigg[
         \int_{-k\tau}^{t \wedge \tau_{n}}
         e^{2\lambda_{1}p
         (s+k\tau)}
         {\big\| X^{-k\tau}_{s}\big\|}^{2p-2}
         \Big(
         \lambda_1
         \| { X^{-k\tau}_{s}\|}^{2}
         +
         \big\langle
         X_{s}^{-k\tau},
         A X_{s}^{-k\tau}
         \big\rangle
         \Big)
         \, \dd s 
         \Bigg] 
         }_{\leq 0} \\
        & \quad 
         +
         2p
         \E 
         \Bigg[
         \int_{-k\tau}
         ^{t \wedge \tau_{n}}
         e^{2\lambda_{1}p
         (s+k\tau)}
         {\big\| X^{-k\tau}_{s}\big\|}^{2p-2}
         \Big(
         \big\langle
         X_{s}^{-k\tau},
         f(s,X_{s}^{-k\tau})
         \big\rangle
         +
         \dfrac{2p-1}{2}
         \big\|
         g(s,X_{s}^{-k\tau}
         )\big\|^2
         \Big)
         \, \dd s 
         \Bigg] \\
        &  \leq
        \E
         [ 
         \|\xi\|^{2p} 
         ] 
         +
         2p 
         \E
         \Bigg[
         \int_{-k\tau}^{t \wedge \tau_{n}}
         e^{2\lambda_{1}p
         (s+k\tau)}
         \Big[
         \alpha_2
         \big\| { X^{-k\tau}_{s}\big\|}^{2p}
         +
         c_0
         \big\| { X^{-k\tau}_{s}\big\|}^{2p-2}
         \Big]
         \, \dd s 
         \Bigg].
\end{split}
\end{equation}
Using the Young inequality 
\begin{equation}
    a^{2p-2}b
    \leq
    \dfrac{p-1}{p}a^{2p}
    +
    \dfrac{1}{p}b^{p},
\end{equation}
{\color{black}for some  positive constant $\epsilon_2 < \lambda_1-\alpha_2$, it can see that}
\begin{equation}
\begin{split}
    c_0 
    \big\| { X^{-k\tau}_{s}\big\|}^{2p-2}
    &=
    (\lambda_1-\alpha_2-\epsilon_2)
    \Big[
    \big\| { X^{-k\tau}_{s}\big\|}^{2p-2}
    \times
    \tfrac{c_o}{\lambda_1-\alpha_2-\epsilon_2}
    \Big] \\
    & \leq
    (\lambda_1-\alpha_2-\epsilon_2)
    \times
    \tfrac{p-1}{p}
    \big\| { X^{-k\tau}_{s}\big\|}^{2p}
    +
    \tfrac{1}{p}
    (\lambda_1-\alpha_2-\epsilon_2)^{1-p}
    c_{0}^{p}\\
    &\leq
    (\lambda_1-\alpha_2-\epsilon_2)
    \big\| { X^{-k\tau}_{s}\big\|}^{2p}
    +
    \tfrac{1}{p}
    (\lambda_1-\alpha_2-\epsilon_2)^{1-p}
    c_{0}^{p}.
\end{split}
\end{equation}
Then one achieves that 
\begin{equation}
\begin{split}
         \E
         \Big[
        e^{2\lambda_{1}p
        (t \wedge \tau_{n}+k\tau)}
        \big\| 
        X^{-k\tau}_{t \wedge \tau_{n}}
        \big\|^{2p}
        \Big] 
        & \leq
        \E
         [ 
         \| \xi\| ^{2p}
         ]
         +
         \E
         \Bigg[
         \int_{-k\tau}^{t \wedge \tau_{n}}
         2 
         (\lambda_1-\epsilon_2) 
         p
         e^{2\lambda_{1}p
         (s+k\tau)}
         \big\| { X^{-k\tau}_{s}\big\|}^{2p}
         \, \dd s 
         \Bigg] \\
         &\quad +
         \E
         \Bigg[
         \int_{-k\tau}^{t \wedge \tau_{n}}
         2
         e^{2\lambda_{1}p
         (s+k\tau)}
         (\lambda_1-\alpha_2-\epsilon_2)^{1-p}
         c_{0}^{p}
         \, \dd s 
         \Bigg] \\
         & \leq
         \E
         [ 
         \| \xi\| ^{2p}
         ]
         +
         \tfrac
         {
         (\lambda_1-\alpha_2-\epsilon_2)^{1-p}
         c_{0}^{p}}
         {\lambda_1 p}
         e^{2\lambda_1 p(t+k\tau)}\\
         & \quad +
         \E\Bigg[
         \int_{-k\tau}^{t \wedge \tau_{n}}
         2p
         (\lambda_1-\epsilon_2) 
         e^{2\lambda_{1}p
         (s+k\tau)}
         \big\| { X^{-k\tau}_{s}\big\|}^{2p}
         \, \dd s 
         \Bigg].
\end{split}
\end{equation}
{\color{black}
By the Gr\"onwall inequality \eqref{eq:Gronwall 1}, we have that
\begin{equation}
\begin{split}
        &
        \E
        \Big[
        e^{2\lambda_{1}p
        (t \wedge \tau_n+k\tau)}
        {\big\| X^{-k\tau}_{t}\big\|}^{2p}
        \Big] \\
        & \quad \leq
        \E
        [\|\xi\|^{2p}]
        +
        \tfrac
         {
         (\lambda_1-\alpha_2-\epsilon_2)^{1-p}
         c_{0}^{p}}
         {\lambda_1 p}
         e^{2\lambda_1 p(t+k\tau)}\\
         & \quad +
         \E\Bigg[
         \int_{-k\tau}^{t \wedge \tau_{n}}
         2 p
         (\lambda_1-\epsilon_2) 
         \Big(
         \| \xi \|^{2p}
         +
        \tfrac
         {
         (\lambda_1-\alpha_2-\epsilon_2)^{1-p}
         c_{0}^{p}}
         {\lambda_1 p}
         e^{2\lambda_1 p(s+k\tau)}
         \Big)
         e^{\int_s^t2p(\lambda_1-\epsilon_2)\,\dd r}
         \, \dd s  
         \Bigg]\\
         & \quad =
         e^{2(\lambda_1-\epsilon_2)p(t+k\tau)}
         \E[\| \xi \|^{2p}]
         +
          \tfrac
         {
         (\lambda_1-\alpha_2-\epsilon_2)^{1-p}
         c_{0}^{p}}
         {\lambda_1 p}
         e^{2\lambda_1 p(t+k\tau)} \\
         & \qquad +
          \tfrac
         {
         (\lambda_1-\epsilon_2)
         (\lambda_1-\alpha_2-\epsilon_2)^{1-p}
         c_{0}^{p}}
         {\lambda_1 \epsilon_2 p}
         e^{2\lambda_1 p(t+k\tau)}
         \big(1
         -
         e^{-2\epsilon_2 p(t+k\tau)}
         \big),
\end{split}
\end{equation}
}
resulting in by letting $n \rightarrow \infty$
\begin{equation}
\begin{split}
        \E
        \Big[
        {\| X^{-k\tau}_{t}\|}^{2p}
        \Big] 
        &\leq
        \E[\| \xi \|^{2p}]
        +
        \tfrac
         {
         (\lambda_1-\alpha_2-\epsilon_2)^{1-p}
         c_{0}^{p}}
         {\lambda_1 p}
         (1+\tfrac{\lambda_1-\epsilon_2}{\epsilon_2}) \\
         & =
         \E[\| \xi \|^{2p}]
        +
        \tfrac
         {
         (\lambda_1-\alpha_2-\epsilon_2)^{1-p}
         c_{0}^{p}}
         {\epsilon_2 p}.
\end{split} 
\end{equation}
The proof is completed.
\end{proof}

We state the following result on the  H\"older continuity of the exact solution of
\eqref{eq_PEM:Problem_SDE} with respect to the norm in 
$L^p(\Omega;\mathbb{R}^d)$.
\begin{lem}
\label{lem_PEM:the_esti_(X(t1)-X(t2)}
    Let
    Assumption \ref{ass_PEM}  
    hold.
    Then there exists a positive constant $C$ which depends on $\gamma,d,A,f,g$ only, such that
    \begin{equation}
\label{eq_PEM:the_esti_(X(t1)-X(t2)}
    \begin{split}
        \big\|
        X_{t_1}^{-k\tau}
        -
        X_{t_2}^{-k\tau}
        \big\|
        _{L^p(\Omega;\mathbb{R}^d)}
        & \leq
        C
        \Big(
          1+\sup_{k\in \mathbb{N}} \sup_{t \geq -k\tau}
          \big\| 
          X_{t}^{-k\tau} 
         \big\|
         ^\gamma_{L^{p\gamma}(\Omega;\mathbb{\R}^d)}
        \Big)
        |t_{2}-t_{1}| \\
        &\quad+
        C
       \Big(
          1+\sup_{k\in \mathbb{N}} \sup_{t\geq -k\tau}
          \big\| 
          X_{t}^{-k\tau} 
         \big\|
         ^{\gamma}
         _{L^{p\gamma}(\Omega;\mathbb{R}^d)}
        \Big)
        |t_{2}-t_{1}|
        ^{\frac{1}{2}} 
        \\
    \end{split}
\end{equation}
    for all $t_1,t_2 \geq -k\tau$ and
    $p \in 
    \Big[2,\frac{2p_1}{\gamma}\Big)$,
    where $X_{t}^{-k\tau}$ denotes the exact solution to the SDE
\eqref{eq_PEM:Problem_SDE}.
\end{lem}

\begin{proof}
[Proof of Lemma \ref{lem_PEM:the_esti_(X(t1)-X(t2)}]
    Without loss of generality we set $t_1 \leq t_2$ and get
    \begin{equation}
    \label{eq_PEM:X_t1-X_t2}
    \begin{split}
        &
        \big\|
        X_{t_1}^{-k\tau}
        -
        X_{t_2}^{-k\tau}
        \big\|
        _{L^p(\Omega;
        \mathbb{R}^d)}
        \\
        &  =
        \bigg\|
        \int_{t_1}^{t_2}
        \big(
        A X_{r}^{k\tau}
        +f(r,X_{r}^{-k\tau}
        \big) 
        \,\dd r
        +
        \int_{t_1}^{t_2} 
        g(r,X_{r}^{-k\tau})
        \,\dd W_r
        \bigg\|
         _{L^p(\Omega;\mathbb{R}^d)} \\
         & \leq
         \bigg\|
         \int_{t_1}^{t_2}
         \big(
           A
           X_{r}^{-k\tau}
           +f(r,X_{r}^{-k\tau})
         \big)
         \,\dd r
         \bigg\|
         _{L^p(\Omega;\mathbb{R}^d)}
          +
          \bigg\|
         \int_{t_1}^{t_2} 
         g(r,X_{r}^{-k\tau})
          \,\dd W_r
         \bigg\|
         _{L^p(\Omega;\mathbb{R}^d)}.
         \end{split}
    \end{equation}
    For the first term, 
    it follows from the H\"older inequality and \eqref{eq_PEM:the_esti_f(t,x)},
    one  can obtain
    \begin{equation}
        \begin{split}
          &
          \bigg\|
          \int_{t_1}^{t_2}
            A X_{r}^{-k\tau} 
            + f(r,X_{r}^{-k\tau})
          \,\dd r
         \bigg\|
         _{L^p(\Omega;\R^d)} \\
         & \quad \leq
         \int_{t_1}^{t_2}
         \big\|
           A X_{r}^{-k\tau}
         \big\|
         _{L^p(\Omega;\R^d)}
         \, \dd r
         +
         \int_{t_1}^{t_2}
        \big\|
         f(r,X_{r}^{-k\tau})
        \big\|
         _{L^p(\Omega;\R^d)}
         \, \dd r \\
         & \quad \leq
         C
         \Big(
               1+\sup_{k \in \N} 
               \sup_{t \geq -k\tau}
               \big\| 
               X_{t}^{-k\tau} 
               \big\|
               ^{\gamma}_{L^{p\gamma}(\Omega;\R^d)}
        \Big)
        |t_2 - t_1|,
        \end{split}
    \end{equation}
    {\color{black}where we make use of Lemma \ref{lem_PEM:the_pth_of_exact_solution} to get the last line}.
    Applying the Burkholder-Davis-Gundy inequality to the last term of \eqref{eq_PEM:X_t1-X_t2} 
    and \eqref{eq:the_esti_g(t,x)} to show
    \begin{equation}
    \begin{split}
        \bigg\|
        \int_{t_1}^{t_2} 
        g(r,X_{r}^{-k\tau})
        \,\dd W_r
        \bigg\|
        _{L^p(\Omega;\mathbb{R}^d)} 
        &\leq
        C 
        \bigg(
        \int_{t_1}^{t_2}
        \big\| 
        g(r,X_{r}^{-k\tau}) 
        \big\|
        ^{2}
        _{L^p(\Omega;\mathbb{R}^d)}
        \, \dd r
        \bigg)^{\frac{1}{2}}\\
        &\leq
        C 
        \Big(
          1+\sup_{k\in \mathbb{N}} \sup_{t\geq -k\tau}
          \big\| 
          X_{t}^{-k\tau} 
         \big\|
         ^{\gamma}
         _{L^{p\gamma}(\Omega;\mathbb{R}^d)}
        \Big)
        |t_2 - t_1 |^{\frac{1}{2}}.
    \end{split}
    \end{equation}
    This completes proof.
\end{proof}

\section{Numerical Approximation of Random Periodic Solution}
\label{sec:PEM}
This section is devoted to the introduction of the projected Euler method for approximating the solution of \eqref{eq_PEM:Problem_SDE} on an infinite horizon.
To do this, consider an equidistant partition
$\mathcal{T}^{h}:={jh,j \in \mathbb{Z}}$, 
such that $h \in(0,1)$.
In addition, for
$x \in \mathbb{R}^{d}$, we define the following function 
\begin{equation} 
\label{eq_PEM:projecion_funcaion}
    \left\{
    \begin{aligned}
    \Phi(x)
    &:=
    \min
    \Big\{
    1,
    h^{-\frac{1}{2\gamma}}
    \|x\|^{-1}
    \Big\}
    x,
    \quad x \neq 0,
    \\
    \Phi(x)
    &=0,
    \quad x = 0,
    \end{aligned}
    \right.
\end{equation}
{\color{black} $\gamma$ is determined in Assumption \ref{ass_PEM} }.

Next we propose our 
explicit numerical method to approximate the exact solution of the 
SDEs \eqref{eq_PEM:Problem_SDE}
starting at $-k\tau$,
\begin{equation}
\label{eq:the_projected_euler_method}
\begin{split}
   {\tilde{X}}
   ^{-k\tau}_{-k\tau+(j+1)h}
   &=
   \Phi
   \big(
   {\tilde{X}}
   ^{-k\tau}_{-k\tau+jh}
   \big)
   +A h
   \Phi
   \big(
   {\tilde{X}}^{-k\tau}_{-k\tau+jh}
   \big)
   + hf
   \Big( 
   jh, 
   \Phi
   \big(
   {\tilde{X}}^{-k\tau}_{-k\tau+jh}
   \big)
   \Big) \\
   &\quad +
   g \Big(
   jh,
   \Phi
   \big(
   {\tilde{X}}^{-k\tau}_{-k\tau+jh}
   \big)
   \Big)
   \Delta W_{-k\tau+jh},
\end{split}
\end{equation}
for all $ j \in \mathbb{N}$, 
where 
$\Delta W_{-k\tau+jh}:=
W_{-k\tau+(j+1)h}-W_{-k\tau+jh}$,
and the initial value 
$\tilde{X}^{-k\tau}_{-k\tau}=\xi$.
Because of the periodicity of $f$ and $g$, we have that 
$f\big(-k\tau+jh,{\tilde{X}}^{-k\tau}_{-k\tau+jh}\big) =f\big(jh,{\tilde{X}}^{-k\tau}_{-k\tau+jh}\big)$, 
$g\big(-k\tau+jh,{\tilde{X}}^{-k\tau}_{-k\tau+jh}\big) =g\big(jh,{\tilde{X}}^{-k\tau}_{-k\tau+jh}\big)$.

Before proceeding further, 
we collect some preliminary estimates,
which have been established in
\cite[Lemma 4.2]{pang2024linear}.
\begin{lem}
 \label{lem_PEM:Phi(x)}
Let Assumption \ref{ass_PEM}
hold,
then for any 
$x \in \mathbb{R}^{d}$,
the following estimates
\begin{equation}
    \|
    \Phi (x)
    \|
    \leq
    h^{-\frac{1}{2\gamma}}, \\
    \quad
    \|
    f(t,\Phi (x))
    \|
    \leq
    L_{1}h^{-\frac{1}{2}},
\end{equation}
hold true,
{\color{black}
where
$L_{1}:=2C_2$, and
$C_2$ is from \eqref{eq_PEM:the_esti_f(t,x)}}.
Moreover, for any
$x,y \in \mathbb{R}^{d}$,
the following estimates hold true
\begin{align}
    \label{eq_PEM:the_esti_Phi(x)-Phi(y)}
    \|
    \Phi(x)-
    \Phi(y)
    \|
    & \leq
    \| x-y \|,
    \\  \label{eq_PEM:the_esti_f(t,Phi(x))-f(t,Phi(y))}
    \|
    f(t,\Phi(x))
    -
    f(t,\Phi(y))
    \|
    & \leq
    L_{2} 
    h^{
    -\frac{\gamma-1}{2\gamma}}
    \|x-y\|,
\end{align}
{\color{black}
where 
$L_{2}:
=3C_1$, and $C_1$ is from \eqref{eq_PEM:the_esti_f(t,x)-f(t,y)}.}
Especially, for $y=0$, we have
for $x \in \R^d$
\begin{equation}
    \label{eq_PEM:Phi(x)}
    \|
    \Phi(x)
    \|
    \leq
    \|x\|.
\end{equation}
\end{lem}

Before proceeding, we present some useful lemma, which has been given in \cite[Lemma 5.7]{pang2023projected}.
\begin{lem}
\label{lem_PEM:x-Phi(x)}
    Recall $\gamma$
    given in Assumption
    \ref{ass_PEM}.
    For all $ x \in \R^d$,
    and let $\Phi(x)$
    be defined as \eqref{eq_PEM:projecion_funcaion}.
    Then when $\gamma >1$,
    one have
    \begin{equation}
        \|
        x
        -
        \Phi(x)
        \|
        \leq
        Ch^2
        \|
        x
        \|
        ^{4\gamma+1}.
    \end{equation}
    In addition,
    when $\gamma+1$,
    $x-\Phi(x)=0$
    for all $ x \in \R^d$.
\end{lem}
The next lemma establishes the following bounds of  time continuity in projection environment.
\begin{lem}
\label{lem_PMM:the_esti_(X(t1)-PhiX(t2)}
    Let
    Assumption \ref{ass_PEM}  
    hold,
    $X_{t}^{-k\tau}$ denotes the exact solution to the SDE
\eqref{eq_PEM:Problem_SDE}.
    Then there exists a positive constant $C$ which depends on $\gamma,d,A,f,g$ only,
    for all $|t_1-t_2|\leq h$ and
    $p \in 
    \Big[2,\frac{2p^*}{\gamma}\Big)$,
    such that
    \begin{equation}
    \label{eq_PMM:the_esti_X(t1)-PhiX(t2)}
        \big\|
        X_{t_1}^{-k\tau}
        -
        \Phi(X_{t_2}^{-k\tau})
        \big\|
        _{L^p(\Omega;\mathbb{R}^d)}
         \leq
        C h^{\frac{1}{2}}
        \Big(
          1+\sup_{k\in \mathbb{N}} \sup_{t \geq -k\tau}
          \big\| 
          X_{t}^{-k\tau} 
         \big\|
         ^{4\gamma+1}
         _{L^{p(4\gamma+1)}(\Omega;\mathbb{\R}^d)}
        \Big).
\end{equation}
\end{lem}
\begin{proof}
[Proof of Lemma \ref{lem_PMM:the_esti_(X(t1)-PhiX(t2)}]
Using triangle inequality leads to
        \begin{equation}
                \big\|
                X_{t_1}^{-k\tau}
                -
                \Phi(X_{t_2}^{-k\tau})
                \big\|
                _{L^p(\Omega;
                \mathbb{R}^d)} \leq
                \big\|
                X_{t_1}^{-k\tau}
                -
                X_{t_2}^{-k\tau}
                \big\|
                _{L^p(\Omega;
                \mathbb{R}^d)}
                +
                \big\|
                X_{t_2}^{-k\tau}
                -
                \Phi(X_{t_2}^{-k\tau})
                \big\|
                _{L^p(\Omega;
                \mathbb{R}^d)}.
        \end{equation}
        Based on Lemma \ref{lem_PEM:the_esti_(X(t1)-X(t2)},
        for $|t_1-t_2| \leq h \in (0,1)$,
        one can get,
        \begin{equation}
        \begin{split}
        \big\|
        &
        X_{t_1}^{-k\tau}
        -
        X_{t_2}^{-k\tau}
        \big\|
        _{L^p(\Omega;\mathbb{R}^d)} 
        \\
        & \leq
        C
        \Big(
        1+\sup_{k\in \mathbb{N}} \sup_{t \geq -k\tau}
        \big\| 
        X_{t}^{-k\tau} 
         \big\|
         ^\gamma_{L^{p\gamma}(\Omega;\mathbb{\R}^d)}
        \Big)
        |t_{2}-t_{1}| 
        +
        C
       \Big(
          1+\sup_{k\in \mathbb{N}} \sup_{t\geq -k\tau}
          \big\| 
          X_{t}^{-k\tau} 
         \big\|
         ^{\gamma}
         _{L^{p\gamma}(\Omega;\mathbb{R}^d)}
        \Big)
        |t_{2}-t_{1}|^{\frac{1}{2}}
        \\
        & \leq
        Ch^{\frac{1}{2}}
        \Big(
          1+\sup_{k\in \mathbb{N}} \sup_{t\geq -k\tau}
          \big\| 
          X_{t}^{-k\tau} 
         \big\|
         ^{\gamma}
         _{L^{p\gamma}(\Omega;\mathbb{R}^d)}
        \Big).
    \end{split}
    \end{equation}
    For the second term,
    applying Lemma \ref{lem_PEM:x-Phi(x)} leads to
    \begin{equation}
        \begin{split}
        \big\|
        X_{t_2}^{-k\tau}
        -
        \Phi(X_{t_2}^{-k\tau})
        \big\|
        _{L^p(\Omega;
        \mathbb{R}^d)}
        & \leq
        Ch^2
        \|
        X_{t_2}^{-k\tau}
        \|
        ^{4\gamma+1}
        _{L^{p(4\gamma+1)}} \\
        & \leq
        Ch^2
        \Big(
          1+\sup_{k\in \mathbb{N}} \sup_{t\geq -k\tau}
          \big\| 
          X_{t}^{-k\tau} 
         \big\|
         ^{4\gamma+1}
         _{L^{p(4\gamma+1)}(\Omega;\mathbb{R}^d)}
        \Big).
        \end{split}
    \end{equation}
    Above all,
    \begin{equation}
        \big\|
        X_{t_1}^{-k\tau}
        -
        \Phi(X_{t_2}^{-k\tau})
        \big\|
        _{L^p(\Omega;\mathbb{R}^d)}
         \leq
        C h^{\frac{1}{2}}
        \Big(
          1+\sup_{k\in \mathbb{N}} \sup_{t \geq -k\tau}
          \big\| 
          X_{t}^{-k\tau} 
         \big\|
         ^{4\gamma+1}
         _{L^{p(4\gamma+1)}(\Omega;\mathbb{\R}^d)}
        \Big).
    \end{equation}
    This completes proof.
\end{proof}
{\color{black}The next lemma shows there is a uniform bound for the second moment of the numerical solution under necessary assumptions.}
\begin{lem}
    \label{lem:second moment of PEM}
    Let Assumption \ref{ass_PEM}
    hold.
    Then there exists a positive constant
    $C$ 
    such that
    \begin{equation}
   \label{eq:second moment of PEM}
        \sup_{k,j \in \mathbb{N}}
        \mathbb{E}
        \Big[
        \big\|
        {\Tilde{X}}
        _{-k\tau+(j+1)h}
        ^{-k\tau}
        \big\|
        ^{2}
        \Big]
        \leq
        C
        \E
        \Big[
        1+
        \|
        \xi
        \|
        ^{2}
        \Big],
   \end{equation}
    where
    $
    \Big\{
    {\Tilde{X}}
    _{-k\tau+(j+1)h}
    ^{-k\tau}
    \Big\}_{k,j \in \N}$
    is given by \eqref{eq:the_projected_euler_method}.
\end{lem}

\begin{proof}
    [Proof of Lemma \ref{lem:second moment of PEM}]
    From the explicit numerical scheme \eqref{eq:the_projected_euler_method},
    we have
    \begin{equation}
    \begin{split}
        &
        \big\|
        {\Tilde{X}}
        _{-k\tau+(j+1)h}
        ^{-k\tau}
        \big\|
        ^{2} \\ 
        & \quad=
        \Big \|
        \Phi
        \big(
        {\tilde{X}}
        ^{-k\tau}_{-k\tau+jh}
       \big)
        +
        A h\Phi
        \big(
        {\tilde{X}}^{-k\tau}_{-k\tau+jh}
        \big)
        + hf
        \Big( 
        jh, 
        \Phi
        \big(
        {\tilde{X}}^{-k\tau}_{-k\tau+jh}
        \big)
        \Big) 
        +
        g
        \Big( 
        jh, 
        \Phi
        \big(
        {\tilde{X}}^{-k\tau}_{-k\tau+jh}
        \big)
        \Big)
        \Delta W_{-k\tau+jh}
        \Big \|
        ^{2}\\
        & \quad=
        \Big \|
         \Phi
        \big(
        {\tilde{X}}
        ^{-k\tau}_{-k\tau+jh}
        \big)
        \Big \|
        ^{2}
          + h^{2}
        \Big \|
        A \Phi
        \big(
        {\tilde{X}}
        ^{-k\tau}_{-k\tau+jh}
        \big)
        \Big\|
        ^{2}
          + h^{2}
        \Big\|
        f
        \Big(
        jh,
        \Phi
        \big(
        {\tilde{X}}
        ^{-k\tau}_{-k\tau+jh}
        \big)
        \Big)
        \Big\|
        ^{2} \\
        & \qquad +
        \Big\|
        g
        \Big(
        jh,
        \Phi
        \big(
        {\tilde{X}}
        ^{-k\tau}_{-k\tau+jh}
        \big)
        \Big)
        \Delta W_{-k\tau+jh}
        \Big\|^2
          + 2h
        \Big \langle
        \Phi
        \big(
        {\tilde{X}}
        ^{-k\tau}_{-k\tau+jh}
        \big),
        A
        \Phi
        \big(
        {\tilde{X}}
        ^{-k\tau}_{-k\tau+jh}
        \big)
        \Big \rangle \\
        & \qquad +
        2h
        \Big \langle
        \Phi
        \big(
        {\tilde{X}}
        ^{-k\tau}_{-k\tau+jh}
        \big),
        f
        \Big(
        jh,
        \Phi
        \big(
        {\tilde{X}}
        ^{-k\tau}_{-k\tau+jh}
        \big)
        \Big)
        \Big \rangle 
        +
        2h^{2}
        \Big \langle
        A
        \Phi
        \big(
        {\tilde{X}}
        ^{-k\tau}_{-k\tau+jh}
        \big),
        f
        \Big(
        jh,
        \Phi
        \big(
        {\tilde{X}}
        ^{-k\tau}_{-k\tau+jh}
        \big)
        \Big)
        \Big \rangle \\
        &\qquad
        +2
        \Big \langle
        \Phi
        \big(
        {\tilde{X}}
        ^{-k\tau}_{-k\tau+jh}
        \big)
        + A h\Phi
        \big(
        {\tilde{X}}^{-k\tau}_{-k\tau+jh}
        \big) 
        + hf
        \big( 
        jh, 
        \Phi
        \big(
        {\tilde{X}}^{-k\tau}_{-k\tau+jh}
        \big),
        g
        \Big(
        jh,
        \Phi
        \big(
        {\tilde{X}}
        ^{-k\tau}_{-k\tau+jh}
        \big)
        \Big)
        \Delta W_{-k\tau+jh}
        \Big \rangle. 
    \end{split}
    \end{equation}
    Taking expectations 
    on both sides and 
    noticing $\E[\Delta W_{-k\tau+jh}|\mathcal{F}_{-k\tau+jh}]=0$,
    one can deduce
\begin{equation}
\label{eq:E(Y^2)}
\begin{split}
        \E
        \Big[
        \big\|
        {\Tilde{X}}
        _{-k\tau+(j+1)h}
        ^{-k\tau}
        \big\|
        ^{2}
        \Big]
        &=
        \E
        \Big[
         \Big \|
         \Phi
        \big(
        {\tilde{X}}
        ^{-k\tau}_{-k\tau+jh}
        \big)
        \Big \|
        ^{2}
        \Big]
        +h^{2}
        \E
        \Big[
        \Big \|
        A
        \Phi
        \big(
        {\tilde{X}}
        ^{-k\tau}_{-k\tau+jh}
        \big)
        \Big\|
        ^{2}
        \Big]
        +h^{2}
        \Big[
        \Big\|
        f
        \big(
        jh,
        \Phi
        \big(
        {\tilde{X}}
        ^{-k\tau}_{-k\tau+jh}
        \big)
        \big)
        \Big\|
        ^{2} 
        \Big] \\
        &\quad +
        h
        \E
        \Big[
        \Big\|
        g
        \Big(
        jh,
        \Phi
        \big(
        {\tilde{X}}
        ^{-k\tau}_{-k\tau+jh}
        \big)
        \Big)
        \Big\| 
        ^{2} 
        \Big]
           + 2h
        \E
        \Big[
        \Big \langle
        \Phi
        \big(
        {\tilde{X}}
        ^{-k\tau}_{-k\tau+jh}
        \big),
        A
        \Phi
        \big(
        {\tilde{X}}
        ^{-k\tau}_{-k\tau+jh}
       \big)
        \Big \rangle
         \Big] \\
        & \quad 
        +2h
        \E\Big[
        \Big \langle
        \Phi
        \big(
        {\tilde{X}}
        ^{-k\tau}_{-k\tau+jh}
        \big),
        f
        \Big(
        jh,
        \Phi
        \big(
        {\tilde{X}}
        ^{-k\tau}_{-k\tau+jh}
        \big)
        \Big)
        \Big \rangle 
        \Big] \\
        &\quad
        +2h^{2}
        \E\Big[
        \Big \langle
        A
        \Phi
        \big(
        {\tilde{X}}
        ^{-k\tau}_{-k\tau+jh}
        \big),
        f
       \Big(
        jh,
        \Phi
        \big(
        {\tilde{X}}
        ^{-k\tau}_{-k\tau+jh}
        \big)
        \Big)
        \Big \rangle
       \Big].
\end{split}
\end{equation}
Using the Young inequality yields
\begin{equation}
\label{eq:2h^2<Ax,f(x)>}
        2h^{2}
        \Big \langle
        A
        \Phi
        \big(
        {\tilde{X}}
        ^{-k\tau}_{-k\tau+jh}
        \big),
        f
        \Big(
        jh,
        \Phi
        \big(
        {\tilde{X}}
        ^{-k\tau}_{-k\tau+jh}
        \big)
        \Big)
        \Big \rangle
        \leq
        h^{2}
       \Big\|
         A
        \Phi
        \big(
        {\tilde{X}}
        ^{-k\tau}_{-k\tau+jh}
        \big)
        \Big\|^{2}
        +
        h^{2}
        \Big\|
        f
        \Big(
        jh,
        \Phi
        \big(
        {\tilde{X}}
        ^{-k\tau}_{-k\tau+jh}
        \big)
        \Big\|
        ^{2}.
\end{equation}
Combing \eqref{eq:2h^2<Ax,f(x)>} into \eqref{eq:E(Y^2)} and
making use of Assumption \ref{ass_PEM}
\eqref{eq_PEM:coercivity_condition}
and Lemma \ref{lem_PEM:Phi(x)} give
\begin{equation}
    \begin{split}
        \E
        \Big[
        \big\|
        {\Tilde{X}}
        _{-k\tau+(j+1)h}
        ^{-k\tau}
        \big\|
        ^{2}
        \Big]
        &\leq
        \E
        \Big[
         \Big \|
         \Phi
        \big(
        {\tilde{X}}
        ^{-k\tau}_{-k\tau+jh}
        \big)
        \Big \|
        ^{2}
        \Big]
        +
        2 h^{2}
        \E
        \Big[
        \Big \|
         A\Phi
        \big(
        {\tilde{X}}
        ^{-k\tau}_{-k\tau+jh}
        \big)
        \Big \|
        ^{2}
        +
        2h^{2}
        \Big\|
        f
        \Big(
        jh,
        \Phi
        \big(
        {\tilde{X}}
        ^{-k\tau}_{-k\tau+jh}
        \big)
        \Big\|
        ^{2}
        \\
        & \quad
         + 2h
        \E
        \Big[
        \Big \langle
        \Phi
        \big(
        {\tilde{X}}
        ^{-k\tau}_{-k\tau+jh}
        \big),
        A
        \Phi
        \big(
        {\tilde{X}}
        ^{-k\tau}_{-k\tau+jh}
        \big)
        \Big \rangle
        \Big] \\
        &\quad +
        2h
        \E \Big[
        \Big \langle
        \Phi
        \big(
        {\tilde{X}}
        ^{-k\tau}_{-k\tau+jh}
        \big),
        f
        \Big(
        jh,
        \Phi
        \big(
        {\tilde{X}}
        ^{-k\tau}_{-k\tau+jh}
        \big)
        \Big)
        \Big \rangle 
         + \tfrac{1}{2}
         \Big\|
        g
        \Big(
        jh,
        \Phi
        \big(
        {\tilde{X}}
        ^{-k\tau}_{-k\tau+jh}
        \big)
        \Big)
        \Big\| 
        ^{2} 
        \Big] \\
        &\leq
        [1-2(\lambda_1-\alpha_2)h]
        \E
        \Big[
        \|
        \Phi
        \big(
        {\Tilde{X}}
        _{-k\tau+jh}^{-k\tau}
        \big)
        \|
        ^{2}
        \Big]
        +2\lambda_{d}^{2}
        h^{2-\frac{1}{\gamma}}
        +
        2 L_1^2 h
        +2 c_0 h\\
        &\leq
        [1-2(\lambda_1-\alpha_2)h]
        \E
        \Big[
        \|
        {\Tilde{X}}
        _{-k\tau+jh}^{-k\tau}
        \|
        ^{2}
        \Big]
        +2\lambda_{d}^{2}
        h^{2-\frac{1}{\gamma}}
        +
        2 L_1^2 h
        +2 c_0 h.
    \end{split}
\end{equation}
{\color{black}Noting that $\gamma \in \Big[1,\frac{p_1+1}{2}\Big)$,
$2\lambda_{d}^{2}
        h^{2-\frac{1}{\gamma}}
        <2\lambda_{d}^{2}h$.}
Combining with Assumption \ref{ass_PEM}
with
$\lambda_1 > \alpha_2$
for some positive $\Tilde{C}:=
C(\lambda_d, L_1, c_0)$,   
such that
\begin{equation}
    \begin{split}
        \E
        \Big[
        \big\|
        {\Tilde{X}}
        _{-k\tau+(j+1)h}
        ^{-k\tau}
        \big\|
        ^{2}
        \Big]
        &\leq
        [1-2(\lambda_{1}-\alpha_{2})h]
        \E
        \Big[
        \big\|
        {\Tilde{X}}
        _{-k\tau+jh}
        ^{-k\tau}
        \big\|
        ^{2}
        \Big] 
        +\Tilde{C}h\\
        &\leq
        [1-2(\lambda_{1}-\alpha_{2})h]^{j+1}
        \E
        \Big[
        \big\|
        {\Tilde{X}}
        _{-k\tau}
        ^{-k\tau}
        \big\|
        ^{2}
        \Big] 
        +
        \sum_{i=0}^{j}
        [1-2(\lambda_{1}-\alpha_{2})h]^{i}
        \Tilde{C}h\\
        &=
        [1-2(\lambda_{1}-\alpha_{2})h]^{j+1}
        \E
        \Big[
        \|
        \xi
        \|
        ^{2}
        \Big] 
        +
        \tfrac{1-[1-2(\lambda_{1}-\alpha_{2})h]^{j+1}}{2(\lambda_{1}-\alpha_{2})h}
        \Tilde{C}h\\
        &\leq
        C
        \E
        \Big[
        1+
        \|
        \xi
        \|
        ^{2}
        \Big].
    \end{split}
\end{equation}
Then the assertion follows.
\end{proof}

The following lemma indicates that any two numerical solutions starting from different initial conditions can be arbitrarily close after sufficiently many iterations.

\begin{lem}
    \label{lem_PEM:the_second_moment_|x-y|}
    Let
    Assumption \ref{ass_PEM} 
    hold
    {\color{black} and recall $L_2$ defined in Lemma \ref{lem_PEM:Phi(x)}}.
    Let
    ${ \tilde{X}^{-k\tau}_{-k\tau+jh}}$ and ${ \tilde{Y}^{-k\tau}_{-k\tau+jh}}$ 
    be two solutions of the 
    projected Euler scheme
   \eqref{eq:the_projected_euler_method}
    with initial values $\xi$
    and $\eta$
    satisfying condition (iv) in Assumption \ref{ass_PEM}, respectively.
    Then it holds that
    \begin{equation}
        \label{eq:the_seond_moment_|x-y|}
        \E
        \Big[
        \big\|
        \tilde{X}^{-k\tau}_{-k\tau+(j+1)h}
        -
        \tilde{Y}^{-k\tau}_{-k\tau+(j+1)h}
        \big\|
        ^{2}
        \Big]
        \leq
        e^{-(\lambda_1-\alpha_1)
        (j+1)h}
        \E
        [
        \| \xi-\eta \|^{2}  ].
    \end{equation}
    {
    \color{black}
    where $h$ is the timestep satisfying
    \begin{equation}
    h 
    \in
    \Bigg(
    0,
    \min
    \Bigg\{
    \frac{(\lambda_1-\alpha_1)^{\frac{p_1+1}{2}}}{(\lambda_{d}+L_2)^{p_1+1}}
    ,
    \frac{1}{\lambda_1-\alpha_1}
    ,
    1
    \Bigg\}
    \Bigg).
    \end{equation}
    }
\end{lem}

\begin{proof}
[Proof of Lemma \ref{lem_PEM:the_second_moment_|x-y|}.]
Subtracting \eqref{eq:the_projected_euler_method} yields
\begin{equation}
\begin{split}
    \tilde{X}^{-k\tau}_{-k\tau+(j+1)h}
    -
    \tilde{Y}^{-k\tau}_{-k\tau+(j+1)h}
    &=
    \Phi
    (\tilde{X}^{-k\tau}_{-k\tau+jh})
    -
    \Phi
    (\tilde{Y}^{-k\tau}_{-k\tau+jh})
    +A h
    \big(
    \Phi
    (\tilde{X}^{-k\tau}_{-k\tau+jh})
    -
    \Phi
    (\tilde{Y}^{-k\tau}_{-k\tau+jh})
    \big)\\
    &\quad+
    h\Big(
    f
    \Big(jh,
    \Phi
    (\tilde{X}^{-k\tau}_{-k\tau+jh}
    \big)
    \Big)
    -
    f
    \Big(jh,
    \Phi
    (\tilde{Y}^{-k\tau}_{-k\tau+jh}
    \big)
    \Big)
    \Big) \\
    &\quad +
    \Big(
    g
    \Big(jh,
    \Phi
    (\tilde{X}^{-k\tau}_{-k\tau+jh}
    \big)
    \Big)
    -
    g
    \Big(jh,
    \Phi
    (\tilde{Y}^{-k\tau}_{-k\tau+jh}
    \big)
    \Big)
    \Big)
    \Delta W_{-k\tau+jh}.
\end{split}
\end{equation}
Shortly, we denote
\begin{align}
\label{eq:x-y}
    \tilde{\zeta}_{j}
    &:=
    \tilde{X}^{-k\tau}_{-k\tau+jh}
    -
    \tilde{Y}^{-k\tau}_{-k\tau+jh},
    \\
    \label{eq:Phi(x)-Phi(y)}
    \Delta 
    \tilde{\Phi}_{j}
    &:=
    \Phi
    (\tilde{X}^{-k\tau}_{-k\tau+jh})
    -
    \Phi
    (\tilde{Y}^{-k\tau}_{-k\tau+jh}),
    \\
    \label{eq:f(t,Phi(x))-f(t,Phi(y))}
    \Delta 
    \tilde{f}_{j} 
    &:=
    f
    \Big(jh,
    \Phi
    (\tilde{X}^{-k\tau}_{-k\tau+jh}
    \big)
    \Big)
    -
    f
    \Big(jh,
    \Phi
    (\tilde{Y}^{-k\tau}_{-k\tau+jh}
    \big)
    \Big),
    \\
    \label{eq:g(t,Phi(x))-g(t,Phi(y))}
    \Delta 
    \tilde{g}_{j} 
    &:=
    g
    \Big(jh,
    \Phi
    (\tilde{X}^{-k\tau}_{-k\tau+jh}
    \big)
    \Big)
    -
    g
    \Big(jh,
    \Phi
    (\tilde{Y}^{-k\tau}_{-k\tau+jh}
    \big)
    \Big).
\end{align}
   With  
\eqref{eq:x-y} to
\eqref{eq:g(t,Phi(x))-g(t,Phi(y))},
it is not hard to show that
\begin{equation}
    \tilde{\zeta}_{j+1}
    =
    \Delta 
    \tilde{\Phi}_{j}
    +A h
    \Delta 
    \tilde{\Phi}_{j}
    +h
    \Delta 
    \tilde{f}_{j}
    +
    \Delta 
    \tilde{g}_{j}
    \Delta W_{-k\tau+jh}.
\end{equation}
Taking the expectation of the second moment on both sides gives
\begin{equation}
\begin{split}
    \E
    \big[
    \|
    \tilde{\zeta}_{j+1}
    \|
    ^{2}
    \big]
    &=
    \E
    \big[
    \big \|
    \Delta 
    \tilde{\Phi}_{j}
    +A h
    \Delta 
    \tilde{\Phi}_{j}
    +h
    \Delta 
    \tilde{f}_{j}
    +
    \Delta 
    \tilde{g}_{j}
    \Delta W_{-k\tau+jh}
    \big \|
    ^{2}
    \big]\\
    &=
    \E
    \big[
    \big\|
    \Delta 
    \tilde{\Phi}_{j}
    \big\|^{2}
    \big]
    +
    h^{2}
    \E
    \big[
    \big\|
    A 
    \Delta 
    \tilde{\Phi}_{j}
    \big\|^{2}
    \big]
    +h^{2}
    \E
    \big[
    \big\|
    \Delta 
    \tilde{f}_{j}
     \big \|^{2}
    \big]
    +
    h
    \E
    \big[
    \|
    \Delta 
    \tilde{g}_{j}
    \|^2
    \big]
    \\
    & \quad+2h
    \E
    \big[
    \langle
    \Delta 
   \tilde{\Phi}_{j},
     A 
    \Delta 
    \tilde{\Phi}_{j}
    \rangle
    \big]
    +2h
    \E
    \big[
    \langle
    \Delta 
    \tilde{\Phi}_{j},
    \Delta 
    \tilde{f}_{j}
    \rangle
    \big]
    +2h^{2}
    \E\big[
    \langle
    A \Delta 
    \tilde{\Phi}_{j},
    \Delta 
    \tilde{f}_{j}
    \rangle
    \big].
\end{split}
\end{equation}
{\color{black}
Using the Cauchy-Schwarz inequality 
$\|\langle a,b\rangle \| \leq \|a\|\|b\|$  leads to
\begin{equation}
    2h^{2}
    \E \big[
    \langle
    A \Delta 
    \tilde{\Phi}_{j},
    \Delta 
    \tilde{f}_{j}
    \rangle
    \big]
    \leq
    2h^{2}
    \E \Big[
    \big\|
    A \Delta 
    \tilde{\Phi}_{j}
    \big\|
    \cdot
    \big\|
    \Delta
    \tilde{f}_{j}
    \big\|
    \Big].
\end{equation}
}
{\color{black}Regarding the terms $\|\Delta \tilde{f}_{j}\|$,
we use Lemma \ref{lem_PEM:Phi(x)} \eqref{eq_PEM:the_esti_f(t,Phi(x))-f(t,Phi(y))} to estimate,}
and recalling Assumption \ref{ass_PEM} and Lemma \ref{lem_PEM:Phi(x)}, one can obtain that
{\color{black}
\begin{equation}\label{eq_PEM:E[zeta_j+1]}
\begin{split}
    \E \big[
    \|
    \tilde{\zeta}_{j+1}
    \|^2
    \big]
    & \leq
    {\color{black}
    \E[
    \|
     \Delta 
    \tilde{\Phi}_{j}
    \|^{2}
    ]}
    +
    2h
    \Big\{
    \E
    \big[
    \langle
    \Delta 
    \tilde{\Phi}_{j},
     A 
    \Delta 
    \tilde{\Phi}_{j}
    \rangle
    \big]
    +
    \E
    \big[
    \langle
    \Delta 
    \tilde{\Phi}_{j},
    \Delta 
    \tilde{f}_{j}
    \rangle
    \big] 
    +
    \tfrac{2p_1-1}{2}
    \E
    \|
    \Delta 
    \tilde{g}_{j}
    \|^2
    \Big\} \\
    & \quad +
    \lambda_{d}^{2}h^{2}
    \E 
    \big[
    \|
     \Delta 
    \tilde{\Phi}_{j}
    \|^{2}
    \big]
    +
    2\lambda_{d}L_{2}
    h^{1+\frac{\gamma+1}{2\gamma}}
    \E 
    \big[
    \|
     \Delta 
    \tilde{\Phi}_{j}
    \|^{2}
    \big] 
    +
    L_{2}^{2}
    h^{1+\frac{1}{\gamma}}
    \E \big[
    \|
     \Delta 
    \tilde{\Phi}_{j}
    \|^{2}
    \big] \\
    &\leq
     \Big(1-2(\lambda_{1}-\alpha_{1})h\Big)
    \E
    \big[
    \|
     \Delta 
    \tilde{\zeta}_{j}
    \|^{2}
    \big]
    +
   \Big(
   \lambda_{d}^{2}h
    + 2\lambda_{d}L_{2}h^{\frac{\gamma+1}{2 \gamma}}
    +
    L_{2}^{2}
    h^{\frac{1}{\gamma}}
    \Big)
    h
    \E \big[
    \|
     \Delta 
    \tilde{\zeta}_{j}
    \|^{2}
    \big]. 
\end{split}
\end{equation}
}
{\color{black}
According to $\gamma \in \Big[1,\frac{p_1+1}{2}\Big)$, one can get
\begin{equation}\label{eq_PEM:h^2/p_1+1}
    \lambda_{d}^{2}h
    + 2\lambda_{d}L_{2}h^{\frac{\gamma+1}{2 \gamma}}
    +
    L_{2}^{2}
    h^{\frac{1}{\gamma}}
    \leq
    (\lambda_{d}+L_2)^2 h^{\frac{2}{p_1 +1}}.
\end{equation}
Here we select an appropriate $h$ such that
\begin{equation}\label{eq_PEM:h<lambda-alpha}
    (\lambda_{d}+L_2)^2 h^{\frac{2}{p_1+1}}
    \leq
    \lambda_1-\alpha_1,
\end{equation}
which leads to
\begin{equation}
    h 
    \in
    \Bigg(
    0,
    \min
    \Bigg\{
    \frac{(\lambda_1-\alpha_1)^{\frac{p_1+1}{2}}}{(\lambda_{d}+L_2)^{p_1+1}}
    ,
    \frac{1}{\lambda_1-\alpha_1}
    ,
    1
    \Bigg\}
    \Bigg).
\end{equation}
Combining \eqref{eq_PEM:h^2/p_1+1} and \eqref{eq_PEM:h<lambda-alpha}
into \eqref{eq_PEM:E[zeta_j+1]},
we can have
\begin{equation}
    \E \big[
    \|
    \tilde{\zeta}_{j+1}
    \|^2
    \big]
    \leq
    \Big(1-(\lambda_{1}-\alpha_{1})h\Big)
    \E
    \big[
    \|
     \Delta 
    \tilde{\zeta}_{j}
    \|^2
    \big].
\end{equation}
As a result,
\begin{equation}
    \E \big[
    \|
    \tilde{\zeta}_{j+1}
    \|
    \big]
    \leq
    \Big(1-(\lambda_{1}-\alpha_{1})h
    \Big)
    \E \big[
    \|
    \tilde{\zeta}_{j}
    \|
    \big]
    \leq
    e^{-(\lambda_{1}-\alpha_{1})(j+1)h}
    \E
    \big[
    \| \xi-\eta \|^{2}
    \big].
\end{equation}
}
Thus we complete the proof.
\end{proof}
{\color{black}
Under the framework of Theorem 3.4 in \cite{feng2017numerical}, we can derive the existence and uniqueness of random periodic solution to the projected Euler method
\eqref{eq:the_projected_euler_method}.}

\begin{thm}
\label{thm:PEM random period solution}
   Let Assumption \ref{ass_PEM}  hold. 
   {\color{black}
   For $h 
    \in
    \Bigg(
    0,
    \min
    \Bigg\{
    \frac{(\lambda_1-\alpha_1)^{\frac{p_1+1}{2}}}{(\lambda_{d}+L_2)^{p_1+1}}
    ,
    \frac{1}{\lambda_1-\alpha_1}
    ,
    1
    \Bigg\}
    \Bigg)$}, 
   the projected Euler method \eqref{eq:the_projected_euler_method} admits a random period solution $\tilde{X}^{*} \in L^{2}(\Omega)$ such that
\begin{equation}
  \lim_{k \rightarrow 
  \infty }
   \E\Big[
   \big\|
   \tilde{X}_{-k\tau+jh}^{-k\tau}(\xi)
   -\tilde{X}^{*}_t
   \big\|^{2}
   \big]
   =0.
\end{equation}
\end{thm}
Note that the Wiener shift $\theta$ appears only in the proof the existence and uniqueness of random periodic solutions of numerical methods. We omit the corresponding proof as it follows the one in 
\cite{feng2017numerical}. In the other words, the error analysis is conducted in a common sense and does not involve the effect of Wiener shift. 

\section{Mean square convergence order of Projected Euler Method}
\label{sec:strong_rate_of_PEM}

We consider the difference between the exact solution and the numerical solution and give a comprehensive error analysis with convergence rate.

The exact solution at time $-k\tau+(j+1)h$ is as follows
\begin{equation}
\label{eq:the_exact_X_j+1}
\begin{split}
        X_{-k\tau+(j+1)h}^{-k\tau} 
        &=
        \Phi
        (X_{-k\tau+jh}^{-k\tau})
        +Ah
        \Phi
        (X_{-k\tau+jh}^{-k\tau})
        +
        hf
        \big(
        jh,
         \Phi(X_{-k\tau+jh}^{-k\tau})
         \big) \\
         &\quad +
          g
          \big(
           jh,
           \Phi(X_{-k\tau+jh}^{-k\tau})
           \big)
           \Delta 
           W_{-k\tau+jh}
           +
           \mathcal{R}
           _{-k\tau+(j+1)h},
\end{split}
\end{equation} 
where,
\begin{equation}
\label{eq:the_def_R_j+1}
    \begin{split}
        \mathcal{R}_{-k\tau+(j+1)h}
        &=
        \int_{-k\tau+jh}^{-k\tau+(j+1)h}
          A
          \big(
          X_{s}^{-k\tau}
          -
          \Phi
          (X_{-k\tau+jh}^{-k\tau})
        \big)
        \, \dd s\\ 
        &\quad +
        \int_{-k\tau+jh}^{-k\tau+(j+1)h}
          f
        \big(
          s,
         X_{s}^{-k\tau}
          \big)
          -
          f
         \big(
          jh,
         \Phi
         (X_{-k\tau+jh}^{-k\tau})
          \big)
        \, \dd s \\
        &\quad +
        \int_{-k\tau+jh}^{-k\tau+(j+1)h}
          g
         \big(
          s,
         X_{s}^{-k\tau}
          \big)
          -
          g
         \big(
          jh,
         \Phi
         (X_{-k\tau+jh}^{-k\tau})
          \big)
        \, \dd W_{s} \\
        & \quad +
          X_{-k\tau+jh}^{-k\tau}
          -
          \Phi
          (
          X_{-k\tau+jh}
          ^{-k\tau}
          ).
    \end{split}
\end{equation}

\subsection{Convergence rates for SDEs with multiplicative noise}
The following lemma provides uniform bounded estimates for the second moment of $\mathcal{R}_{-k\tau+(j+1)h}$
and its conditional expectation $\E[\mathcal{R}_{-k\tau+(j+1)h}|\mathcal{F}_{-k\tau+jh}]$.
\begin{lem}
\label{lem_PEM:the_esti_E(R^2)}
    Let Assumption \ref{ass_PEM} hold.
    Then for $k,j \in \N$,
    there exists some positive constant $C$,
    independent of $k,j$
    {\color{black} and h,}
    such that
    \begin{equation}
        \|
        \mathcal{R}_{-k\tau+(j+1)h}
        \|_{L^2(\Omega;\R^d)}
        \leq
        Ch,
        \qquad
        \big\|
        \E[
        \mathcal{R}_{-k\tau+(j+1)h}
        | \mathcal{F}
        _{-k\tau+jh}
        ]\big\|_{L^2(\Omega;\R^d)}
        \leq
        Ch^{\frac{3}{2}}.
    \end{equation}
\end{lem}
\begin{proof}
[Proof of Lemma \ref{lem_PEM:the_esti_E(R^2)}]
    Recalling the definition of 
    $\mathcal{R}_{-k\tau+(j+1)h}$ given by
    \eqref{eq:the_def_R_j+1}
    and using an triangle inequality yield
    \begin{equation}
    \label{eq:R=I_1+I_2+I_3}
    \begin{split}
        \|
        \mathcal{R}_{-k\tau+(j+1)h}
        \|
        _{L^2(\Omega;\R^d)}
        & \leq
        \bigg\|
        \int_{-k\tau+jh}^{-k\tau+(j+1)h}
        A\big(
          X_{s}^{-k\tau}
          -
          \Phi
          (X_{-k\tau+jh}^{-k\tau})
        \big)
        \, \dd s 
        \bigg\|
        _{L^2(\Omega;\R^d)}\\
        &\quad +
        \bigg\|
        \int_{-k\tau+jh}^{-k\tau+(j+1)h}
          f
         \big(
          s,
         X_{s}^{-k\tau}
          \big)
          -
          f
         \big(
          jh,
         \Phi
         (X_{-k\tau+jh}^{-k\tau})
          \big)
        \,\dd s 
        \bigg\|
        _{L^2(\Omega;\R^d)}\\
        &\quad +
        \bigg\|
        \int_{-k\tau+jh}^{-k\tau+(j+1)h}
          g
         \big(
          s,
         X_{s}^{-k\tau}
          \big)
          -
          g
         \big(
          jh,
         \Phi
         (X_{-k\tau+jh}^{-k\tau})
          \big)
        \, \dd W_{s}
        \bigg\|
        _{L^2(\Omega;\R^d)} \\
        & \quad +
        \|
          X_{-k\tau+jh}^{-k\tau}
          -
          \Phi
          (
          X_{-k\tau+jh}
          ^{-k\tau}
          )
          \|_{L^2(\Omega;\R^d)}
        \\
        &
        := 
        \I_1
        +\I_2
        +\I_3
        +\I_4. 
    \end{split}
\end{equation}
    For the term $\I_1$,
    if follows the H\"older inequality
    and 
    \eqref{eq_PMM:the_esti_X(t1)-PhiX(t2)}
    to give
    \begin{equation}
    \label{eq_PEM:the_esti_I_1}
        \begin{split}
            \I_1
            & \leq
            \int_{-k\tau+jh}
            ^{-k\tau+(j+1)h}
            \big\|
            A
            \big(
            X_{s}^{-k\tau}
            -
            \Phi
            (X_{-k\tau+jh}^{-k\tau})
            \big)
            \big\|
            _{L^2(\Omega;\R^d)}
            \, \dd s \\
            & \leq
            Ch^{\tfrac{3}{2}}
            \Big(
             1+\sup_{k\in \mathbb{N}} \sup_{t\geq -k\tau}
              \big\| 
              X_{t}^{-k\tau} 
             \big\|
             ^{4\gamma+1}
             _{L^{8\gamma+2}(\Omega;\mathbb{R}^d)}
             \Big).
        \end{split}
    \end{equation}
    To estimate $\I_2$,
    using the H\"older
    inequality yields
    \begin{equation}
    \begin{split}
            \I_2
            & =
            \bigg\|
            \int_{-k\tau+jh}
            ^{-k\tau+(j+1)h}
            f
            \big(s,
            X_{s}^{-k\tau}
            \big)
            -
            f
            \big(
            jh,
            \Phi
            (X_{-k\tau+jh}^{-k\tau})
            \big)
            \, \dd s
            \bigg\|
            _{L^2(\Omega;\R^d)}
            \\
            & \leq
            \int
            _{-k\tau+jh}
            ^{-k\tau+(j+1)h}
            \big\|
            f
            \big(
            s,
            X_{s}^{-k\tau}
            \big)
            -
            f
            \big(
            jh,
            X_{s}^{-k\tau}
            \big)
            \big\|
            _{L^2(\Omega;\R^d)}
            \, \dd s \\
            & \quad +
             \int
             _{-k\tau+jh}
             ^{-k\tau+(j+1)h}
            \big\|
            f
            \big(
            jh,
            X_{s}^{-k\tau}
            \big)
            -
            f
            \big(
            jh,
            \Phi
            (
            X_{-k\tau+jh}^{-k\tau}
            )
            \big)
            \big\|
            _{L^2(\Omega;\R^d)}
            \, \dd s \\
            & =:
            \I_{21}+\I_{22}.
    \end{split}
    \end{equation}
    For the term $\I_{21}$,
    applying \eqref{eq_PEM:the_esti_f(t,x)-f(s,y)} leads to
    \begin{equation}
    \begin{split}
        \I_{21}
        & \leq
        \int
        _{-k\tau+jh}
        ^{-k\tau+(j+1)h}
        C
        \Big(
        1
        +
        \|
        X_{s}^{-k\tau}
        \|
        ^{\gamma}
        _{L^{2\gamma}(\Omega;\mathbb{\R}^d)}
        \Big)
        |
        s
        -
        jh
        |
        \, \dd s \\
        & \leq
        Ch^2
        \Big(
          1+\sup_{k\in \mathbb{N}} \sup_{t \geq -k\tau}
          \big\| 
          X_{t}^{-k\tau} 
         \big\|
         ^\gamma
         _{L^{2\gamma}(\Omega;\mathbb{\R}^d)}
        \Big).
    \end{split}
    \end{equation}
    For the term $\I_{22}$,
    it follows from \eqref{eq_PEM:the_esti_f(t,x)-f(t,y)} and \eqref{eq_PEM:Phi(x)}
    \begin{equation}
    \begin{split}
        \I_{22}
        &\leq
        \int_{-k\tau+jh}
        ^{-k\tau+(j+1)h}
        \bigg\|
            C
            \Big(
            1+
            \big\|
            X_{s}^{-k\tau}
            \big\|^{\gamma-1}
            +
            \big\|
            \Phi
            (X_{-k\tau+jh}^{-k\tau})
            \big\|^{\gamma-1}
            \Big)
            \big\|
            X_{s}^{-k\tau}
            -\Phi
            (X_{-k\tau+jh}^{-k\tau})
            \big\| 
            \bigg\|
            _{L^2(\Omega;\R^d)}
            \, \dd s \\
            &\leq
            \int_{-k\tau+jh}
            ^{-k\tau+(j+1)h}
            \bigg\|
            C
            \Big(
            1+
            \big\|
            X_{s}^{-k\tau}
            \big\|^{\gamma-1}
            +
            \big\|
            X_{-k\tau+jh}
            ^{-k\tau}
            \big\|^{\gamma-1}
            \Big)
            \big\|
            X_{s}^{-k\tau}
            -\Phi
            (X_{-k\tau+jh}^{-k\tau})
            \big\| 
            \bigg\|
            _{L^2(\Omega;\R^d)}
            \, \dd s.
    \end{split}
    \end{equation}
    Using the H\"older inequality    
    $$\|v^{\gamma-1}u\|_{L^2(\Omega;\mathbb{R}^d)} 
    \leq
    \|v\|^{\gamma-1}_{L^{2\rho_1(\gamma-1)}(\Omega;\mathbb{R}^d)}
    {\color{black}\times}
    \|u\|_{L^{2\rho_2}(\Omega;\mathbb{R}^d)}
    , $$
    {\color{black}
    for $\frac{1}{\rho_1}+\frac{1}{\rho_2}=1$ with exponents }
    $\rho_1 :=
    \frac{5\gamma}{\gamma-1}$ and
    $\rho_2:=
    \frac{5\gamma}{4\gamma+1}$ 
    yields that
    \begin{equation}
            \I_{22}
            \leq
            \int_{-k\tau+jh}
            ^{-k\tau+(j+1)h}
            C
            \Big(
            1+
            \sup_{k \in \N}
            \sup_{t \geq -k\tau}
            \big\|
            X_{t}^{-k\tau}
            \big\|
            ^{\gamma-1}
            _{L^{2\rho_1(\gamma-1)}(\Omega;\R^d)}
            \Big)
            \big\|
            X_{s}^{-k\tau}
            -
            \Phi
            (
            X_{-k\tau+jh}
            ^{-k\tau}
            )
            \big\|
            _{L^{2\rho_2}(\Omega;\R^d)}.
    \end{equation}
    Moreover, through 
\eqref{eq_PMM:the_esti_X(t1)-PhiX(t2)}
   with $p=2\rho_2$ we have that
   \begin{equation}
           \| 
           X_{s}^{-k\tau}
           -
           \Phi
           (
           X_{-k\tau+jh}^{-k\tau} 
           )
           \|_{L^{2\rho_2}(\Omega;\mathbb{\R}^d)}
           \leq
           C
           h^{\frac{1}{2}}
          \Big(
           1+
           \sup_{k\in \N} 
           \sup_{t\geq -k\tau}
            \| X_{t}^{-k\tau}\| ^{4\gamma+1}
            _{L^{10\gamma}(\Omega;\R^d)}
            \Big)
   \end{equation}
    Note that 
    $2\rho_1(\gamma-1)
    =10\gamma$. Altogether, it follows that for $s \in [-k\tau+jh,-k\tau+(j+1)h]$
    \begin{equation}
            \I_{22}
            \leq
            Ch^{\frac{3}{2}}
            \Big(
            1+
            \sup_{k\in \N} \sup_{t\geq -k\tau}
            \| X_{t}^{-k\tau} \|
            ^{5\gamma}
            _{L^{10\gamma}(\Omega;
            \R^d)}\Big).
    \end{equation}
    Above all,
    \begin{equation}
    \label{eq_PEM:the_esti_I_2}
        \I_2
        \leq
         C
         h^{\frac{3}{2}}
            \Big(
            1+
           \sup_{k\in \N} 
           \sup_{t
           \geq -k\tau}
            \| X_{t}^{-k\tau} \|
            ^{5\gamma}
            _{L^{10\gamma}(\Omega;\R^d)}
            \Big).
    \end{equation}
    For the term $\I_3$,
    in view of the It\^o isomery,
    we get
    \begin{equation}
        \I_3
        =
        \Bigg(
        \int_{-k\tau+jh}^{-k\tau+(j+1)h}
        \big\|
          g
         \big(
          s,
         X_{s}^{-k\tau}
          \big)
          -
          g
         \big(
          jh,
         \Phi
         (X_{-k\tau+jh}^{-k\tau})
          \big)
        \big\|
        ^2
        _{L^{2}(\Omega;\R^d)}
        \, \dd s
        \Bigg)^{\frac{1}{2}}.
    \end{equation}
    Similarly, one also
    obtains
    \begin{equation}
        \I_3
        \leq
        C
        h
            \Big(
            1+
            \sup_{k \in \N} \sup_{t\geq -k\tau}
            \| X_{t}^{-k\tau} \|
            ^{5\gamma}
            _{L^{10\gamma}(\Omega;\R^d)}
            \Big).
    \end{equation}
    With regard to $\mathbb{I}_4$,
    using Lemma \ref{lem_PEM:x-Phi(x)} leads to
    \begin{equation}
       \mathbb{I}_4
       \leq
       Ch^2
       \|
       X_{-k\tau+jh}
       ^{-k\tau}
       \|
       ^{4\gamma+1}
       _{L^{8\gamma+2}(\Omega;\R^d)}.
   \end{equation}
    Putting all the above estimates together we derive from \eqref{eq:R=I_1+I_2+I_3}
    that
    \begin{equation}
        \|
        \mathcal{R}_{-k\tau+(j+1)h}
        \|
        _{L^2(\Omega;\R^d)}
        \leq
        C
        h
            \Big(
            1+
            \sup_{k \in \N} \sup_{t\geq -k\tau}
            \| X_{t}^{-k\tau} \|
            ^{5\gamma}
            _{L^{10\gamma}(\Omega;\R^d)}
            \Big).
    \end{equation}
    Note that 
    $\E
        \Big[
         \int_{-k\tau+jh}^{-k\tau+(j+1)h}
          g
         \big(
          s,
         X_{s}^{-k\tau}
          \big)
          -
          g
         \big(
          jh,
         \Phi
         (X_{-k\tau+jh}^{-k\tau})
          \big)
        \, \dd W_s 
        \Big|
        \mathcal{F}_{-k\tau+jh}
        \Big]
        =0$.
    Using the Jensen inequality for conditional
    expectation to get
    \begin{equation}
        \begin{split}
        &
        \big\|
        \E
        [
        \mathcal{R}_{-k\tau+(j+1)h}
        | \mathcal{F}
        _{-k\tau+jh}
        ]
        \big\|_{L^2(\Omega;\R^d)} \\
        & \quad \leq 
        \bigg\|
        \E
        \bigg[
        \int_{-k\tau+jh}^{-k\tau+(j+1)h}
        A
        \big(
        X_{s}^{-k\tau}
        -
        \Phi
        (X_{-k\tau+jh}
        ^{-k\tau})
        \big)
        \, \dd s 
        \Big|
        \mathcal{F}_{-k\tau+jh}
        \bigg]
        \bigg\|_{L^2(\Omega;\R^d)}\\
        & \qquad  +
        \bigg\|
        \E
        \bigg[
         \int_{-k\tau+jh}^{-k\tau+(j+1)h}
          f
         \big(
          s,
         X_{s}^{-k\tau}
          \big)
          -
          f
         \big(
          jh,
         \Phi
         (X_{-k\tau+jh}^{-k\tau})
          \big)
        \, \dd s 
        \Big|
        \mathcal{F}_{-k\tau+jh}
        \bigg]
        \bigg\|_{L^2(\Omega;\R^d)}\\
        & \qquad  +
        \|
        X_{-k\tau+jh}
        ^{-k\tau}
        -\Phi
        (
        X_{-k\tau+jh}
        ^{-k\tau}
        )
        \|_{L^2(\Omega;\R^d)}
        \\
        & \quad \leq 
        \bigg\|
        \int_{-k\tau+jh}^{-k\tau+(j+1)h}
        A
        \big(
        X_{s}^{-k\tau}
        -
        \Phi
        (X_{-k\tau+jh}
        ^{-k\tau})
        \big)
        \, \dd s 
        \bigg\|_{L^2(\Omega;\R^d)}\\
        & \qquad  +
        \bigg\|
         \int_{-k\tau+jh}^{-k\tau+(j+1)h}
          f
         \big(
          s,
         X_{s}^{-k\tau}
          \big)
          -
          f
         \big(
          jh,
         \Phi
         (X_{-k\tau+jh}^{-k\tau})
          \big)
        \, \dd s 
        \bigg\|_{L^2(\Omega;\R^d)} \\
        & \qquad  +
        \|
        X_{-k\tau+jh}
        ^{-k\tau}
        -\Phi
        (
        X_{-k\tau+jh}
        ^{-k\tau}
        )
        \|_{L^2(\Omega;\R^d)}.
        \end{split}
    \end{equation}
    Recalling \eqref{eq_PEM:the_esti_I_1},\eqref{eq_PEM:the_esti_I_2}
    and lemma \ref{lem_PEM:x-Phi(x)}
    it immediately follows that
    \begin{equation}
        \big\|
        \E
        [
        \mathcal{R}_{-k\tau+(j+1)h}
        | \mathcal{F}
        _{-k\tau+jh}
        ]\big\|_{L^2(\Omega;\R^d)}
        \leq
        Ch^{\frac{3}{2}}
        \Big(
            1+
            \sup_{k \in \N} \sup_{t\geq -k\tau}
            \| X_{t}^{-k\tau} \|
            ^{5\gamma}
            _{L^{10\gamma}(\Omega;\R^d)
            }
            \Big).
    \end{equation}
\end{proof}
We are now ready to give the main result of this section that reveals the  convergence of the projected Euler scheme to the SDE \eqref{eq_PEM:Problem_SDE} in the long run.

\begin{thm}
\label{thm_PEM:error analysis}
Let Assumptions
\ref{ass_PEM}  hold
{\color{black} and recall $L_2$ defined in Lemma \ref{lem_PEM:Phi(x)}}.
If  $X_{-k\tau+jh}^{-k\tau}$ and 
$\tilde{X}_{-k\tau+jh}^{-k\tau}$ 
are the exact and the numerical solutions given by \eqref{eq_PEM:Problem_SDE} and \eqref{eq:the_projected_euler_method}, respectively.
{\color{black}For an arbitrary
pair $(\delta_1,\delta_2)$ s.t.
$\delta_1\in(0,\lambda_1-\alpha_1)$ and
$\delta_2>0$,}
then
there exists a positive constant $C$,
{\color{black}independent of $k,j$ and $h$}, such that
\begin{equation}
    \sup_{k,j \in \N}
    \E
    \Big[
    \big\|
    X_{-k\tau+(j+1)h}
    ^{-k\tau}
    -
    \tilde{X}_{-k\tau+(j+1)h}
    ^{-k\tau}
    \big\|^2
    \Big]
    \leq
    Ch,
\end{equation}
{\color{black}
where the timestep $h$ satisfies
   \begin{equation}
        h 
    \in
    \Bigg(
    0,
    \min
    \Bigg\{
    \frac{(\lambda_1-\alpha_1)^{\frac{p_1+1}{2}}}
    {
    (1+\delta_2)
    ^{\frac{p_1+1}{2}}(\lambda_{d}+L_2)^{p_1+1}}
    ,
    \frac{1}{\lambda_1-\alpha_1-\delta_1}
    ,
    1
    \Bigg\}
    \Bigg).
\end{equation}
}
\end{thm}
\begin{proof}
    [Proof of Theorem \ref{thm_PEM:error analysis}]
    Recalling  \eqref{eq:the_projected_euler_method}
    and   \eqref{eq:the_exact_X_j+1}
    yields
    \begin{equation}
    \label{eq_PEM:X_j+1-Y_j+1}
    \begin{split}
        & X_{-k\tau+(j+1)h}^{-k\tau}
        -
        \tilde{X}_{-k\tau+(j+1)h}^{-k\tau}\\
        & =
        \Phi
        (X_{-k\tau+jh}^{-k\tau})
        -
        \Phi
        (\tilde{X}_{-k\tau+jh}^{-k\tau}) 
         +
        Ah
        \big(
         \Phi
         (X_{-k\tau+jh}^{-k\tau})
         -
         \Phi
         (\tilde{X}_{-k\tau+jh}^{-k\tau})
        \big) \\
        &\quad +
        h
        \Big[
          f
         \big(
          jh,
         \Phi
         (X_{-k\tau+jh}^{-k\tau})
          \big)
          -
          f
         \big(
          jh,
         \Phi
         (\tilde{X}_{-k\tau+jh}^{-k\tau})
          \big)
        \Big] \\
        & \quad +
        \Big[
          g
        \big(
          jh,
         \Phi
         (X_{-k\tau+jh}^{-k\tau})
          \big)
          -
          g
         \big(
          jh,
         \Phi
         (\tilde{X}_{-k\tau+jh}^{-k\tau})
          \big)
        \Big] 
        \Delta W_{-k\tau+jh}
        +
        \mathcal{R}_{-k\tau+(j+1)h}.
    \end{split}
    \end{equation}
    For brevity, we denote
    \begin{equation}
    \label{eq_PEM:the_def_e_j}
    \begin{split}
          e_{-k\tau+jh}
          &:=
           X_{-k\tau+jh}^{-k\tau}
           -
          \tilde{X}_{-k\tau+jh}^{-k\tau},
          \\
          \Delta \Phi_{-k\tau+jh}^{x}
          &:=
          \Phi
          (X_{-k\tau+jh}^{-k\tau})
          -
          \Phi
          (\tilde{X}_{-k\tau+jh}^{-k\tau}), 
          \\
          \Delta \Phi_{{-k\tau+jh}}^{f}
          &:=
           f
         \big(
          jh,
         \Phi
         (X_{-k\tau+jh}^{-k\tau})
          \big)
          -
          f
         \big(
          jh,
         \Phi
         (\tilde{X}_{-k\tau+jh}^{-k\tau})
          \big),
          \\
          \Delta \Phi_{{-k\tau+jh}}^{g}
          &:=
           g
         \big(
          jh,
         \Phi
         (X_{-k\tau+jh}^{-k\tau})
          \big)
          -
          g
         \big(
          jh,
         \Phi
         (\tilde{X}_{-k\tau+jh}^{-k\tau})
          \big),
    \end{split}
    \end{equation}
    we emphasize that 
    $\Delta \Phi_{{-k\tau+jh}}^{x}$,
    $\Delta \Phi_{{-k\tau+jh}}^{f}$,
    and
    $\Delta \Phi_{{-k\tau+jh}}^{g}$
    are 
    $\mathcal{F}_{-k\tau+jh}$-measurable.
    Using \eqref{eq_PEM:the_def_e_j},
    \eqref{eq_PEM:X_j+1-Y_j+1}
     can be rewritten as
    \begin{equation}
        e_{-k\tau+(j+1)h}
        =
        \Delta \Phi_{{-k\tau+jh}}^{x}
        +h 
        \big(
          A \Delta \Phi_{{-k\tau+jh}}^{x}
          +
          \Delta \Phi_{{-k\tau+jh}}^{f}
        \big)
        +
        \Delta \Phi_{{-k\tau+jh}}^{g}
        \Delta W_{-k\tau+jh}
        +
        \mathcal{R}_{-k\tau+(j+1)h}.
    \end{equation}
    This leads to
    \begin{equation}
    \begin{split}
        \|
        e_{-k\tau+(j+1)h}
        \|^2
        & =
        \|
         \Delta \Phi_{{-k\tau+jh}}^{x}
        \|^2
        +
        h^2
        \|
          A \Delta \Phi_{{-k\tau+jh}}^{x}
          +
          \Delta \Phi_{{-k\tau+jh}}^{f}
        \|^2
        \\&\quad+
        \|
        \Delta \Phi_{{-k\tau+jh}}^{g}
        \Delta W_{-k\tau+jh}
        \|^2
        +
        \|
        \mathcal{R}_{-k\tau+(j+1)h}
        \|^2 \\
        &\quad
        +
        2h
        \langle
         \Delta \Phi_{{-k\tau+jh}}^{x},
         A \Delta \Phi_{{-k\tau+jh}}^{x}
          +
          \Delta \Phi_{-k\tau+jh}^{f}
        \rangle 
        + 
        2
        \langle
         \Delta \Phi_{-k\tau+jh}^{x},
         \Delta \Phi_{-k\tau+jh}^{g}
         \Delta W_{-k\tau+jh}
        \rangle \\
        & \quad+
        2
        \langle
         \Delta \Phi_{-k\tau+jh}^{x},
         \mathcal{R}_{-k\tau+(j+1)h}
        \rangle 
        +2h
        \langle
         A\Delta \Phi_{-k\tau+jh}^{x}
         +
         \Delta \Phi_{-k\tau+jh}^{f},
         \Delta \Phi_{-k\tau+jh}^{g}
         \Delta W_{-k\tau+jh}
         \rangle \\
         & \quad +2h
         \langle
         A\Delta \Phi_{-k\tau+jh}^{x}
         +
         \Delta \Phi_{-k\tau+jh}^{f},
         \mathcal{R}_{-k\tau+(j+1)h}
        \rangle 
        +
        2
        \langle
         \Delta \Phi_{-k\tau+jh}^{g}
         \Delta W_{-k\tau+jh},
         \mathcal{R}_{-k\tau+(j+1)h}
        \rangle .
    \end{split}
    \end{equation}
    Taking expectations on both sides gives
    \begin{equation} 
    \label{eq_PEM:E|e_j+1|^2}
        \begin{split}
            &\E
            [
            \|
            e_{-k\tau+(j+1)h}
            \|^2
            ] \\
            & \quad=
            \E
            [
            \|
            \Delta \Phi_{-k\tau+jh}^{x}
            \|^2
            ]
            +
            h^2
            \E
            [
            \|
            A \Delta \Phi_{-k\tau+jh}^{x}
            +
            \Delta \Phi_{-k\tau+jh}^{f}
            \|^2
            ]
            +
            h
            \E
            [
            \|
            \Delta \Phi_{-k\tau+jh}^{g}
            \|^2
            ] 
            +
            \E
            [
            \|
            \mathcal{R}_{-k\tau+(j+1)h}
            \|^2 
            ] \\
            & \qquad +
            2
            \E[
            \langle
            \Delta \Phi_{-k\tau+jh}^{x},
            h
            (
            A \Delta \Phi_{-k\tau+jh}^{x}
            +
            \Delta \Phi_{-k\tau+jh}^{f}
            )
            \rangle
            ]
            +
            2
            \E
            [
            \langle
            \Delta \Phi_{-k\tau+jh}^{x},
            \mathcal{R}_{-k\tau+(j+1)h}
            \rangle 
            ] \\
            & \qquad
            +2
            \E[
            \langle
            h
            (A\Delta \Phi_{-k\tau+jh}^{x}
            +
            \Delta \Phi_{-k\tau+jh}^{f}),
            \mathcal{R}_{-k\tau+(j+1)h}
            \rangle 
            ]
            +
            2
            \E[
            \langle
            \Delta \Phi_{-k\tau+jh}^{g}
            \Delta W_{-k\tau+jh},
            \mathcal{R}_{-k\tau+(j+1)h}
            \rangle 
            ].
        \end{split}
    \end{equation}
    Due to $\Delta \Phi^x_{-k\tau+jh}$
    is $\mathcal{F}_{-k\tau+jh}$-measurable,
    applying the Cauchy-Schwartz inequality
    $2ab \leq \delta_1 ha^2+\tfrac{1}{\delta_1 h}b^2$
     with 
     {\color{black}$0<\delta_1 <\lambda_1-\alpha_1$} for arbitrary positive $h$,
    we deduce
    \begin{equation}
    \label{eq_PEM:2E<Phi^x,R>}
    \begin{split}
            2
            \E
            [
            \langle
            \Delta \Phi_{-k\tau+jh}^{x},
            \mathcal{R}_{-k\tau+(j+1)h}
            \rangle
            ]
            &=
            2
            \E
            [
            \E
            \langle
            \Delta \Phi_{-k\tau+jh}^{x},
            \mathcal{R}_{-k\tau+(j+1)h}
            \rangle
            |
            \mathcal{F}_{-k\tau+jh}
            ]
            ]
            \\
            &=
            2
            \E
            [
            \langle
            \Delta \Phi_{-k\tau+jh}^{x},
            \E
            [
            \mathcal{R}_{-k\tau+(j+1)h}
            |
            \mathcal{F}_{-k\tau+jh}
            ]
            \rangle
            ]
            \\
            &\leq
            \delta_1 h
            \E[
            \|
            \Delta \Phi_{-k\tau+jh}^{x}
            \|^2
            ]
            +
            \tfrac{1}{\delta_1 h}
            \E
            [
            \|
            \E
            [
            \mathcal{R}_{-k\tau+(j+1)h}
            |
            \mathcal{F}_{-k\tau+jh}
            ]
            \|^2
            ].
    \end{split}
    \end{equation}
   Regarding the seventh term of \eqref{eq_PEM:E|e_j+1|^2}, 
   {\color{black}
   for a positive $\delta_2$},
   using the Young inequality leads to
    \begin{equation}
    \label{eq_PEM:2E<h(APhi^x+Phi^f,R)>}
    \begin{split}
            &2
            \E[
            \langle
            h
            (A\Delta \Phi_{-k\tau+jh}^{x}
            +
            \Delta \Phi_{-k\tau+jh}^{f}),
            \mathcal{R}_{-k\tau+(j+1)h}
            \rangle 
            ] \\
            & \quad \leq
            \delta_2 h^2
            \E[
            \|
            A\Delta \Phi_{-k\tau+jh}^{x}
            +
            \Delta \Phi_{-k\tau+jh}^{f}
            \|^2
            ]
            +
            \tfrac{1}{\delta_2}
            \E
            [
            \|
            \mathcal{R}_{-k\tau+(j+1)h}
            \|^2
            ].
    \end{split}
    \end{equation}
    Similarly, one also obtains
    \begin{equation}
    \label{eq_PEM:2E<Phi^g,R>}
            2
            \E[
            \langle
            \Delta \Phi_{-k\tau+jh}^{g}
            \Delta W_{-k\tau+jh},
            \mathcal{R}_{-k\tau+(j+1)h}
            \rangle 
            ]
            \leq
            (2{\color{black}p_1}-2)h
            \E
            [
            \|
            \Delta
            \Phi_{-k\tau+jh}^g
            \|^2
            ]
            +
            \tfrac{1}{2{\color{black}p_1}-2}
            \E
            [
            \|
            \mathcal{R}_{-k\tau+(j+1)h}
            \|^2
            ].
    \end{equation}
    Substituting \eqref{eq_PEM:2E<Phi^x,R>},
    \eqref{eq_PEM:2E<h(APhi^x+Phi^f,R)>} and 
    \eqref{eq_PEM:2E<Phi^g,R>} into \eqref{eq_PEM:E|e_j+1|^2} and using
    \eqref{eq_PEM:lambda},
    \eqref{eq:coupled__momotoncity_condition}
    yield
    \begin{equation}
        \begin{split}
           & \E
            [
            \|
            e_{-k\tau+(j+1)h}
            \|^2
            ]\\
            & \leq
            (1+\delta_1 h)
            \E
            [
            \|
            \Delta \Phi_{-k\tau+jh}^{x}
            \|^2
            ]
            +
            (1+\delta_2)
            h^2
            \E
            [
            \|
            A \Delta \Phi_{-k\tau+jh}^{x}
            +
            \Delta \Phi_{-k\tau+jh}^{f}
            \|^2
            ] \\
            & \quad+
            2
            h
            \E[
            \langle
            \Delta \Phi_{-k\tau+jh}^{x},
            A \Delta \Phi_{-k\tau+jh}^{x}
            +
            \Delta \Phi_{-k\tau+jh}^{f}
            \rangle
            ]
            +
            (2{\color{black}p_1}-1)
            h
            \E
            [
            \|
            \Delta \Phi_{-k\tau+jh}^{g}
            \|^2
            ] \\
            & \quad+
            (
            1+\tfrac{1}{\delta_2}
            +\tfrac{1}{2{\color{black}p_1}-2})
            \E
            [
            \|
            \mathcal{R}_{-k\tau+(j+1)h}
            \|^2 
            ] 
            +
            \tfrac{1}{\delta_1 h}
            \E
            [
            \|
            \E[
            \mathcal{R}_{-k\tau+(j+1)h}
            |
            \mathcal{F}_{-k\tau+jh}
            ]
            \|^2
            ] \\
            & \leq
            [1-(2\lambda_1-2\alpha_1-\delta_1 )h]
            \E
            [
            \|
            \Delta \Phi_{-k\tau+jh}^{x}
            \|^2
            ]
            +
            (1+\delta_2)
            h^2
            \E
            [
            |
            A \Delta \Phi_{-k\tau+jh}^{x}
            +
            \Delta \Phi_{-k\tau+jh}^{f}
            \|^2
            ] \\
            & \quad +
            (
            1+\tfrac{1}{\delta_2}
            +\tfrac{1}{2{\color{black}p_1}-2})
            \E
            [
            \|
            \mathcal{R}_{-k\tau+(j+1)h}
            \|^2 
            ] 
            +
            \tfrac{1}{\delta_1 h}
            \E
            [
            \|
            \E[
            \mathcal{R}_{-k\tau+(j+1)h}
            |
            \mathcal{F}_{-k\tau+jh}
            ]
            \|^2
            ].
        \end{split}
    \end{equation}
    {\color{black}
    Note that, applying 
    Assumption \ref{ass_PEM}
    and 
    Lemma \ref{lem_PEM:Phi(x)}
    and
    using the same technique in \eqref{eq_PEM:h^2/p_1+1}
    to get
    \begin{equation}
    \label{eq_PEM:h(Ax+f)}
        \begin{split}
            & h^2
            \|
             A \Delta \Phi_{-k\tau+jh}^{x}
            +
            \Delta \Phi_{-k\tau+jh}^{f}
            \|^2 \\
            & \quad \leq
            \lambda_{d}^2
            h^2
            \|
            \Delta \Phi_{-k\tau+jh}^{x}
            \|^2
            +
            2
            \lambda_{d} L_2
            h^{1+\frac{\gamma+1}{2\gamma}}
            \|
            \Delta \Phi_{-k\tau+jh}^{x}
            \|^2
            +
            L_2^2
            h^{1+\frac{1}{\gamma}}
            \|
            \Delta \Phi_{-k\tau+jh}^{x}
            \|^2 \\
            & \quad \leq
            \Big(
            (
            \lambda_{d} 
            +
            L_2
            )^2 h^{\frac{2}{p_1+1}}
            \Big)h
            \|
            \Delta \Phi_{-k\tau+jh}^{x}
            \|^2.
        \end{split}
    \end{equation}
    Here we select an appropriate $h$
    leads to 
    \begin{equation}
        h 
    \in
    \Bigg(
    0,
    \min
    \Bigg\{
    \frac{(\lambda_1-\alpha_1)^{\frac{p_1+1}{2}}}
    {
    (1+\delta_2)
    ^{\frac{p_1+1}{2}}(\lambda_{d}+L_2)^{p_1+1}}
    ,
    \frac{1}{\lambda_1-\alpha_1-\delta_1}
    ,
    1
    \Bigg\}
    \Bigg),
    \end{equation}
    to ensure
    \begin{equation}
    (1+\delta_2)
    (
    \lambda_{d} 
    +
    L_2
    )^2 h^{\frac{2}{p_1+1}}
    <
    \lambda_1-\alpha_1,
    \quad
    1-\big(\lambda_1-\alpha_1-\delta_1
    \big)h
    >0.
    \end{equation}
    Above all,
    \begin{equation}
    \label{eq:e_j+1}
        \begin{split}
            \E
            [
            \|
            e_{-k\tau+(j+1)h}
            \|^2
            ]
            & \leq
            \Big\{
            1-
            \big(\lambda_1-\alpha_1-\delta_1
            \big)
            h
            \Big\}
            \E
            [
            \|
            e_{-k\tau+jh}
            \|^2
            ]
            +
            (
            1+\tfrac{1}{\delta_2}
            +\tfrac{1}{2p_1-2})
            \E
            [
            \|
            \mathcal{R}_{-k\tau+(j+1)h}
            \|^2 
            ] \\
            & \quad +
            \tfrac{1}{\delta_1 h}
            \E
            [
            \E
            [
            \|
            \mathcal{R}_{-k\tau+(j+1)h}
            |
            \mathcal{F}_{-k\tau+jh}
            ]
            \|^2
            ].
        \end{split}
    \end{equation}
Denote 
    $c:=\lambda_{1}
    -\alpha_{1}
    -\delta_1$,
    and recall that $\E
        [
        \|
        \mathcal{R}_{-k\tau+(j+1)h}
        \|^2
        ]
        \leq
        Ch^2$ and 
        $\E
        [
        \|
        \E[
        \mathcal{R}_{-k\tau+(j+1)h}
        | \mathcal{F}
        _{-k\tau+jh}
        ]
        \|^2
        ]
        \leq
        Ch^3$ from Lemma \ref{lem_PEM:the_esti_E(R^2)}, we have that      
    \begin{equation}
        \begin{split}
            \E
            [
            \|
            e_{-k\tau+(j+1)h}
            \|^2
            ]
            & \leq
            (1-ch)
            \E
            [
            \|
            e_{-k\tau+jh}
            \|^2
            ]+
            Ch^2 \\
            & \leq
            (1-ch)^{j+1}
            \E
            [\|e_{-k\tau}\|^2]
            +
            \sum_{i=0}^j
            (1-ch)^i
            Ch^2\\
            &=
            (1-ch)^{j+1}
            \E
            [\|e_{-k\tau}\|^2]
            +
            \tfrac{1-(1-ch)^{j+1}}{ch}
            Ch^2.
        \end{split}
    \end{equation}
    }
    {\color{black}
    By observing 
    $e_{-k\tau}=X^{-k\tau}_{-k\tau}-\tilde{X}^{-k\tau}_{-k\tau}=0$,
    one can deduce
    \begin{equation}
        \E
        [\|
        e_{-k\tau+(j+1)h}
        \|^2]
        \leq
        Ch,
    \end{equation}
    then the assertion follows.
    }
\end{proof}

\subsection{Convergence rates for SDEs with additive noise}

In the present subsection,
If SDEs \eqref{eq_PEM:Problem_SDE}
driven by additive noise,
taking the form of
\begin{equation}
	\label{eq:add_SDE}
	\left\{
	\begin{aligned}
	\dd X_{t}^{t_{0}} 
         & = 
         \big(
         A X_{t}^{t_{0}}
         +
         f(t, X_{t}^{t_{0}})
         \big)
         \,\dd t 
         + 
         g(t)
         \, \dd W_t,
   \quad t \in (t_{0},T], \\
		 X_{t_{0}}^{t_{0}} & = \xi.
	\end{aligned}\right.
\end{equation}
Now we revisit a necessary assumption in \cite{guo2023order}.
\begin{assumption}
\label{ass:add}
    Suppose the diffusion coefficient functions {\color{black}$g:\R \rightarrow \R$ are continuous and periodic in time with period $\tau >0$,
    i.e.,
    $g(t+\tau)=g(t)$
    for all $t \in \R$}.
    Besides, there exists a constant $c_g >0$
    such that
    $\sup_{t \in [0,\tau)}
    \|g(t)\| \leq c_g$
    and 
    \begin{equation}
        \|
        g(t_1)
        -g(t_2)
        \|
        \leq
        c_g|t_2-t_1|,
        \quad
        \forall
        t_1, t_2 \in 
        [0,\tau).
    \end{equation}
    Moreover,
    assume the drift coefficient functions
    $f:\R \times \R^d
    \rightarrow \R^d$
    are continuously differentiable,
    and there exists a constant 
    {\color{black}$\gamma \in \Big[1,\frac{p_1+1}{2}\Big)$} and $c_f <\lambda_1$
    such
    that
    \begin{align}
       \langle
       x-y,
       f(x)-f(y)
       \rangle 
       &\leq 
       c_f
       \|x-y\|^2,
       \quad
        \forall
        x,y \in \R^d,
        \\
        \big\|
        \big(
        \tfrac{\partial f}{\partial x}
        (t,x)
        -
        \tfrac{\partial f}{\partial y}
        (t,y)
        \big)
        z
        \big\|
        &\leq
        C
        (
        1+
        \|x\|
        +
        \|y\|
        )
        ^{\max\{\gamma-2,0\}}
        \|x-y\|
        \|z\|,
        \quad
        \forall
        x,y,z \in \R^d,
        t \in [0,\tau),
    \end{align}
{\color{black}
where $\tfrac{\partial f}{\partial x}$ denotes the partial derivative of $f$
with respect to the state variable x.}
\end{assumption}

Based on the above assumption,
{\color{black}we can 
improve the estimates in Lemma \ref{lem_PEM:the_esti_E(R^2)}
by the following lemma given in \cite[Theorem 4.6]{guo2023order}.}
The proof of the following lemma is thus omitted.

\begin{lem}
\label{lem_PEM:the_esti_add_E(R^2)}
    Let Assumptions \ref{ass_PEM} 
    and \ref{ass:add} hold.
    Then for $k,j \in \N$,
    there exists some positive constant $C$,
    independent of $k,j$
    {\color{black}and h},
    such that
    \begin{equation}
        \|
        \mathcal{R}_{-k\tau+(j+1)h}
        \|_{L^2(\Omega;\R^d)}
        \leq
        Ch^{\frac{3}{2}},
        \qquad
        \|
        \E[
        \mathcal{R}_{-k\tau+(j+1)h}
        | \mathcal{F}
        _{-k\tau+jh}
        ]
        \|_{L^2(\Omega;\R^d)}
        \leq
        Ch^2.
    \end{equation}
\end{lem}

\begin{thm}
\label{thm_PEM:add_error analysis}
Let Assumptions
\ref{ass_PEM} and \ref{ass:add} hold.
If  $X_{-k\tau+jh}^{-k\tau}$ and 
$\tilde{X}_{-k\tau+jh}^{-k\tau}$ 
are the exact and the numerical solutions given by \eqref{eq:add_SDE} and \eqref{eq:the_projected_euler_method}, respectively. 
{\color{black}For an arbitrary
pair $(\mu_1,\mu_2)$ s.t.
$\mu_1\in(0,\lambda_{1}-c_f)$ and 
$\mu_2>0$,
}
then
there exists a constant $C>0$,
{\color{black}independent of $k,j$ and $h$}, such that
\begin{equation}
    \sup_{k,j \in \N}
    \E
    [
    \|
    X_{-k\tau+jh}
    ^{-k\tau}
    -
    \tilde{X}_{-k\tau+jh}
    ^{-k\tau}
    \|^2
    ]
    \leq
    Ch^2,
\end{equation}
{\color{black}
where the timestep $h$ satisfies
   \begin{equation}
        h 
    \in
    \Bigg(
    0,
    \min
    \Bigg\{
    \frac{(\lambda_1-c_f)^{\frac{p_1+1}{2}}}
    {
    (1+\mu_2)^{\frac{p_1+1}{2}}(\lambda_{d}+L_2)^{p_1+1}}
    ,
    \frac{1}{\lambda_1-c_f-\mu_1}
    ,
    1
    \Bigg\}
    \Bigg).
\end{equation}
}
\end{thm}

\begin{proof}
    [Proof of Theorem \ref{thm_PEM:add_error analysis}]
    {\color{black}
    Repeating \eqref{eq_PEM:X_j+1-Y_j+1} used in Theorem \ref{thm_PEM:error analysis},
    the term $\Delta \Phi^g_{-k\tau+jh}$ disappears,
    one can get
     \begin{equation}
        e_{-k\tau+(j+1)h}
        =
        \Delta \Phi_{{-k\tau+jh}}^{x}
        +h 
        \big(
          A \Delta \Phi_{{-k\tau+jh}}^{x}
          +
          \Delta \Phi_{{-k\tau+jh}}^{f}
        \big)
        +
        \mathcal{R}_{-k\tau+(j+1)h}.
    \end{equation}
    We emphasize that 
    $\Delta \Phi_{{-k\tau+jh}}^{x}$
    {\color{black}and}
    $\Delta \Phi_{{-k\tau+jh}}^{f}$
    are $\mathcal{F}_{-k\tau+jh}$-
    measurable.
    Taking the expectation of the second moment on both sides gives
     \begin{equation} 
        \begin{split}
            &\E
            [
            \|
            e_{-k\tau+(j+1)h}
            \|^2
            ] \\
            & =
            \E
            [
            \|
            \Delta \Phi_{-k\tau+jh}^{x}
            \|^2
            ]
            +
            h^2
            \E
            [
            \|
            A \Delta \Phi_{-k\tau+jh}^{x}
            +
            \Delta \Phi_{-k\tau+jh}^{f}
            \|^2
            ]
            +
            \E
            [
            \|
            \mathcal{R}_{-k\tau+(j+1)h}
            \|^2 
            ] \\
            & \quad +
            2
            \E[
            \langle
            \Delta \Phi_{-k\tau+jh}^{x},
            h
            (
            A \Delta \Phi_{-k\tau+jh}^{x}
            +
            \Delta \Phi_{-k\tau+jh}^{f}
            )
            \rangle
            ]
            +
            2
            \E
            [
            \langle
            \Delta \Phi_{-k\tau+jh}^{x},
            \mathcal{R}_{-k\tau+(j+1)h}
            \rangle 
            ] \\
            & \quad
            +2
            \E[
            \langle
            h
            (A\Delta \Phi_{-k\tau+jh}^{x}
            +
            \Delta \Phi_{-k\tau+jh}^{f}),
            \mathcal{R}_{-k\tau+(j+1)h}
            \rangle 
            ].
        \end{split}
    \end{equation}
    {\color{black}Recall that} $\Delta \Phi^x_{-k\tau+jh}$ is $\mathcal{F}_{-k\tau+jh}$-measurable
    {\color{black}and note that $2ab\leq \mu_1 h a^2+\tfrac{1}{\mu_1 h}b^2$ for arbitrary $h>0$, $\mu_1\in(0,\lambda_{1}-c_f)$},
    we obtain
    \begin{equation}
    \begin{split}
            2
            \E
            [
            \langle
            \Delta \Phi_{-k\tau+jh}^{x},
            \mathcal{R}_{-k\tau+(j+1)h}
            \rangle
            ]
            &=
            2
            \E
            [
            \E
            \langle
            \Delta \Phi_{-k\tau+jh}^{x},
            \mathcal{R}_{-k\tau+(j+1)h}
            \rangle
            |
            \mathcal{F}_{-k\tau+jh}
            ]
            ]
            \\
            &=
            2
            \E
            [
            \langle
            \Delta \Phi_{-k\tau+jh}^{x},
            \E
            [
            \mathcal{R}_{-k\tau+(j+1)h}
            |
            \mathcal{F}_{-k\tau+jh}
            ]
            \rangle
            ]
            \\
            &\leq
            \mu_1 h
            \E[
            \|
            \Delta \Phi_{-k\tau+jh}^{x}
            \|^2
            ]
            +
            \tfrac{1}{\mu_1 h}
            \E
            [
            \|
            \E
            [
            \mathcal{R}_{-k\tau+(j+1)h}
            |
            \mathcal{F}_{-k\tau+jh}
            ]
            \|^2
            ].
    \end{split}
    \end{equation} 
   For a positive $\mu_2$,
   applying the Young inequality yields
    \begin{equation}
    \begin{split}
            &2
            \E[
            \langle
            h
            (A\Delta \Phi_{-k\tau+jh}^{x}
            +
            \Delta \Phi_{-k\tau+jh}^{f}),
            \mathcal{R}_{-k\tau+(j+1)h}
            \rangle 
            ] \\
            & \quad \leq
            \mu_2 h^2
            \E[
            \|
            A\Delta \Phi_{-k\tau+jh}^{x}
            +
            \Delta \Phi_{-k\tau+jh}^{f}
            \|^2
            ]
            +
            \tfrac{1}{\mu_2}
            \E
            [
            \|
            \mathcal{R}_{-k\tau+(j+1)h}
            \|^2
            ].
    \end{split}
    \end{equation}
    {\color{black}
    Due to \eqref{eq_PEM:h(Ax+f)},
    applying Assumptions 
    \ref{ass_PEM}, \ref{ass:add}
    and Lemma \ref{lem_PEM:Phi(x)},
    one can get
     \begin{equation}
        \begin{split}
            \E
            [
            \|
            e_{-k\tau+(j+1)h}
            \|^2
            ]
            & \leq
            \Big \{
            1-
            \Big(
            \lambda_1-c_f-\mu_1
            \Big)h
            \Big\}
            \E
            [
            \|
            e_{-k\tau+jh}
            \|^2
            ]
            +
            (
            1+\tfrac{1}{\mu_2})
            \E
            [
            \|
            \mathcal{R}_{-k\tau+(j+1)h}
            \|^2 
            ] \\
            & \quad +
            \tfrac{1}{\mu_1 h}
            \E
            [
            \E
            \big(
            \|
            \mathcal{R}_{-k\tau+(j+1)h}
            \|
            \mathcal{F}_{-k\tau+jh}
            \big)
            \|^2
            ]. 
        \end{split}
    \end{equation}
Denoting $\tilde{c}:=\lambda_{1}
    -c_f
    -\mu_1$,
    and taking  lemma \ref{lem_PEM:the_esti_add_E(R^2)}, result in 
    \begin{equation}
            \E
            [
            \|
            e_{-k\tau+(j+1)h}
            \|^2
            ]
            \leq
            (1-\tilde{c}h)
            \E
            [
            \|
            e_{-k\tau+jh}
            \|^2
            ]+
            Ch^3,
    \end{equation}
    Now using a similar argument as the proof of
    Theorem \ref{thm_PEM:error analysis}
    we can deduce that}
    \begin{equation}
        \E
        [\|
        e_{-k\tau+(j+1)h}
        \|^2]
        \leq
        Ch^2,
    \end{equation}
which completes the proof.
}
\end{proof}

\begin{cor}
\label{cor:Ch}
Let Assumption \ref{ass_PEM} 
hold,
let $X^{*}_{t}$ be the random periodic solution of SDE \eqref{eq_PEM:Problem_SDE} and $\tilde{X}^{*}_{t}$ be the random periodic solution of the projected Euler numerical approximation.
Then there exists a constant $C$ independent of {\color{black}$t$ and $h$}, such that 
\begin{equation}
  \sup_{t \in \mathcal{T}^{h}} 
  \E([\| X^{*}_{t}-\tilde{X}^{*}_{t}\|^{2}])^{1/2}
  \leq Ch^{\frac{1}{2}},
\end{equation}
If in addition Assumption 
\ref{ass:add} hold,
then there exists $C>0$,
independent of {\color{black}$t$ and $h$},
such that
\begin{equation}
  \sup_{t \in \mathcal{T}^{h}} 
  \E([\| X^{*}_{t}-\tilde{X}^{*}_{t}\|^{2}])^{1/2}
  \leq Ch.
\end{equation}
\end{cor}

\begin{proof}
[Proof of Corollary \ref{cor:Ch}]
Due to
\begin{equation}
\E[
\| X^{*}_{t}-\tilde{X}^{*}_{t}\|^{2}]
\leq
\limsup_{k}
\big[
\E[\| X^{*}_{t}-X^{-k\tau}_{t}\|^{2}]
+\E[\| X^{-k\tau}_{t}-\tilde{X}^{-k\tau}_{t}\|^{2}]
+\E[\|\tilde{X}^{-k\tau}_{t}-\tilde{X}^{*}_{t}\|^{2}]\big],
\end{equation}
thus the conclusion can be obtained by Theorem \ref{thm:unique_random_periodic_solution}, Theorem \ref{thm:PEM random period solution},
Theorem \ref{thm_PEM:error analysis}
and Theorem \ref{thm_PEM:add_error analysis}.
\end{proof}

\section{Numerical experiments}
\label{sec:numerical results}
Some numerical experiments will be performed to illustrate the previous theoretical results in this section.
To accomplish this, we consider two examples of SDEs with multiplicative noise and additive noise.

\subsection{Example 1}\label{sec:eg1}
\begin{figure}[H]
   \centering
\includegraphics[scale=0.5]
{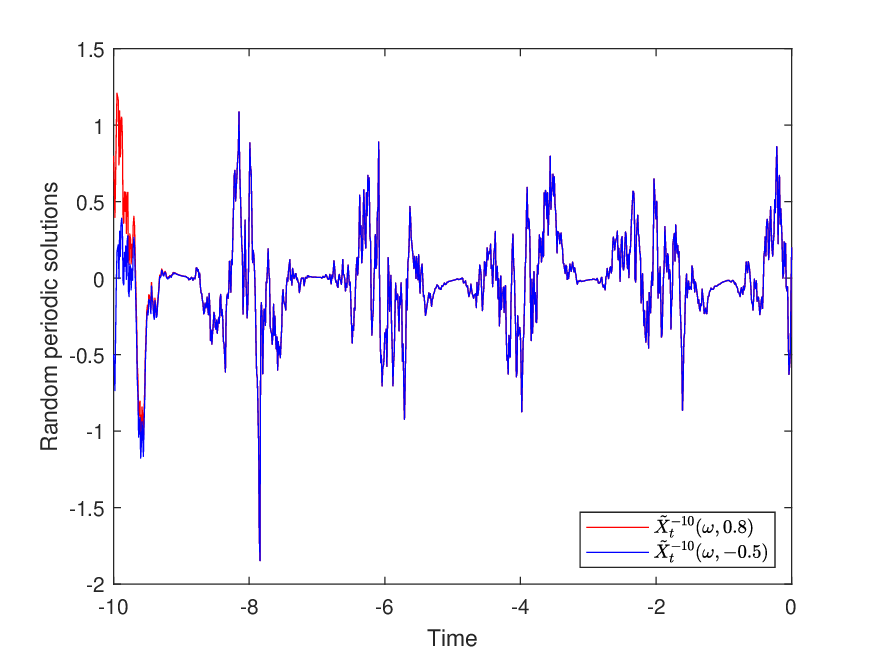}
	\caption[Figure 1.]
	{Two paths generated by projected Euler methods from differential initial conditions.}
 \label{fig1}
\end{figure}

In the first example, we test the performance of the projected Euler method \eqref{eq:the_projected_euler_method} with multiplicative noise as follows:
\begin{equation}
\label{eq_PEM:example_1}
\dd X^{t_{0}}_{t}
=\Big(-2\pi X^{t_{0}}_{t}
+
X^{t_{0}}_{t}
-(X^{t_{0}}_{t})^3
+\cos(\pi t)\Big)
\, \dd t
+
\Big(1+(X^{t_{0}}_{t})^2+
\cos(\pi t)\Big)
\, \dd W_{t},
X^{t_0}_{t_0}=\xi.
\end{equation}
{\color{black}It's easy to verify that \eqref{eq_PEM:example_1}
satisfies Assumptions \ref{ass_PEM}.}
Building upon this,
Theorem \ref{thm:PEM random period solution} states 
its projected Euler simulation also displays a random periodic path.
To further validate this claim,
we conduct numerical experiment 
where we observe two processes starting from $t_0=-10$ and $T=0$,
with the stepsize of 0.01 and initial values of $0.8$ and $-0.5$.
Figure \ref{fig1} illustrates that two paths converge quickly,
demonstrating that the random periodic solution of the projected Euler methods is independent of the initial values.

\begin{figure}[H]
   \centering
\includegraphics[scale=0.5]
{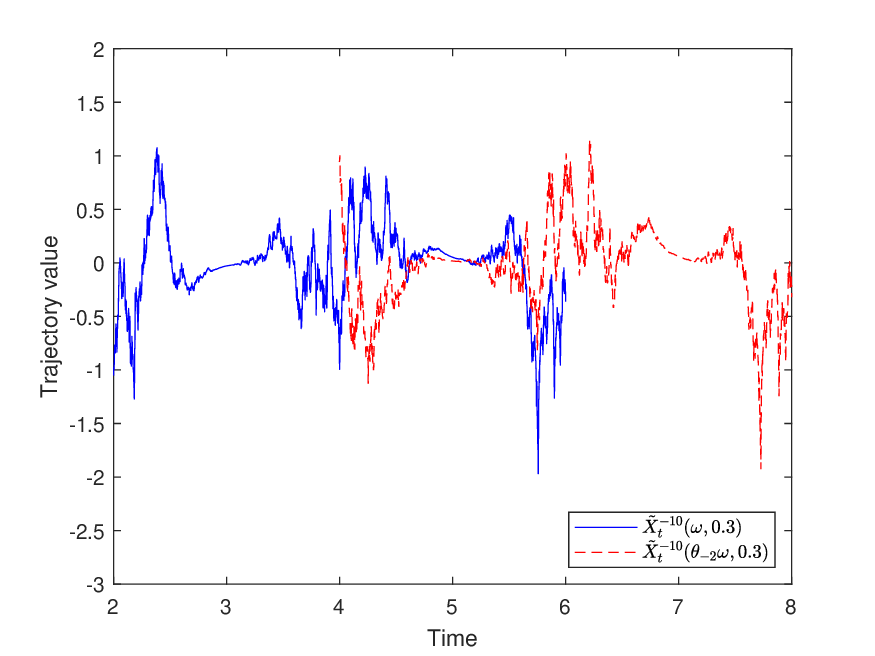}
\caption[Figure 2.]
	{Simulations of the process 
 ${\tilde{X}_{t}^{-10}(\omega,0.3),2 \leq t \leq 6}$
 and
 ${\tilde{X}_{t}^{-10}(\theta_{-2}\omega,0.3),4 \leq t \leq 8}$.}
 \label{fig2}
\end{figure}

Next, we verify the periodicity by examining and dynamics under the same realisation $\omega:\tilde{X}_{t}^{-10}(\omega,0.3)$ over $-10 \leq t \leq 6$
and $\tilde{X}_{t}^{-10}(\theta_{-2}\omega,0.3)$ over $-10 \leq t \leq 8$,
where $0.3$ is the initial condition of both processes.
Due to Theorem \ref{thm:PEM random period solution},
 it is expected that $\tilde{X}^{-10}_{t}(\omega,0.3) \approx \tilde{X}^{*}_{t}(\omega)$ and $\tilde{X}^{-10}_{t}(\theta_{-2}\omega,0.3)\approx \tilde{X}^{*}_{t}(\theta_{-2}\omega)$ after a sufficiently long time, and we may then observe that $\tilde{X}^{-10}_{t-2}(\omega,0.3)\approx \tilde{X}^{-10}_{t}(\theta_{-2}\omega,0.3)$ due to the fact $\tilde{X}^{*}_{t-2}(\omega)=\tilde{X}^{*}_{t}(\theta_{-2}\omega)$ in Definition \ref{def:rps}.
 Figure \ref{fig2} demonstrates both process resemble each other with a stable time gap 2, that is,
 $\tilde{X}^{*}_{t-2}(\omega)=\tilde{X}^{*}_{t}(\theta_{-2}\omega)$ over
 $4 \leq t \leq 8 $.

{\color{black}
Theorem \ref{thm_PEM:error analysis} suggests
that the random periodic solution converges to the solution of \eqref{eq_PEM:example_1} with order 0.5 in the mean square sense.
To achieve this,
 a fine stepsize $h_{ecact}=2^{-15}\times 20 $ is chosen to obtain 
a reference solution on the time interval $[-10,10]$.
The reference solution is obtained via the same numerical method with a fine stepsize $h_{exact}=2^{-15} \times 20$.
 We plot mean-square approximation errors 
$e_{h}$ against five different stepsizes $h=2^{-i} \times 20$, $i=8,9,...,12$ on a log-log scale.}

\begin{figure}[H]
	\centering	\includegraphics[scale=0.5]
 {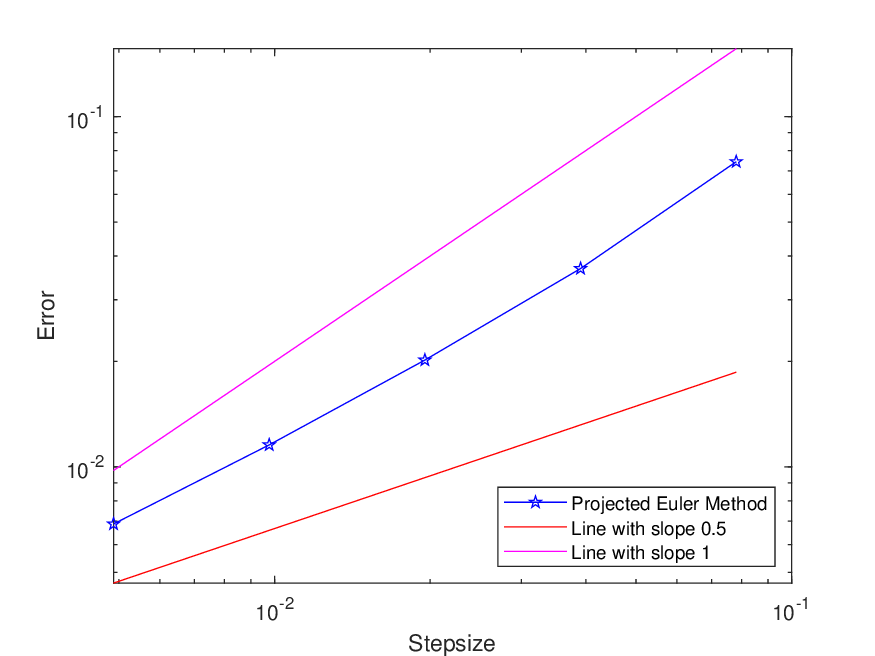}	
 \caption[Figure 3.]
{The mean-square error plot of the projected Euler method \eqref{eq:the_projected_euler_method} for simulating the solution of \eqref{eq_PEM:example_1}.}
\label{fig3}
\end{figure}

{\color{black}
Figure \ref{fig3} clearly demonstrates that the mean-square error is at a slope 
greater than $0.5$, 
but less than $1$.}
Suppose that the approximation error $e_{h}$ obeys a 
power law relation $e_{h}=Ch^{\kappa}$ for 
$C, \kappa >0$, so that $\log e_{h}=\log C+\kappa \log h$.
Then we do a least squares power law fit for 
$\kappa $ and get the value
0.8544 for the rate $\kappa $  with  residual of 0.1112, which is beyond the theoretical
order of convergence 
in Theorem \ref{thm_PEM:error analysis}.

\subsection{Example 2}
\begin{figure}[H]
	\centering
\includegraphics[scale=0.5]{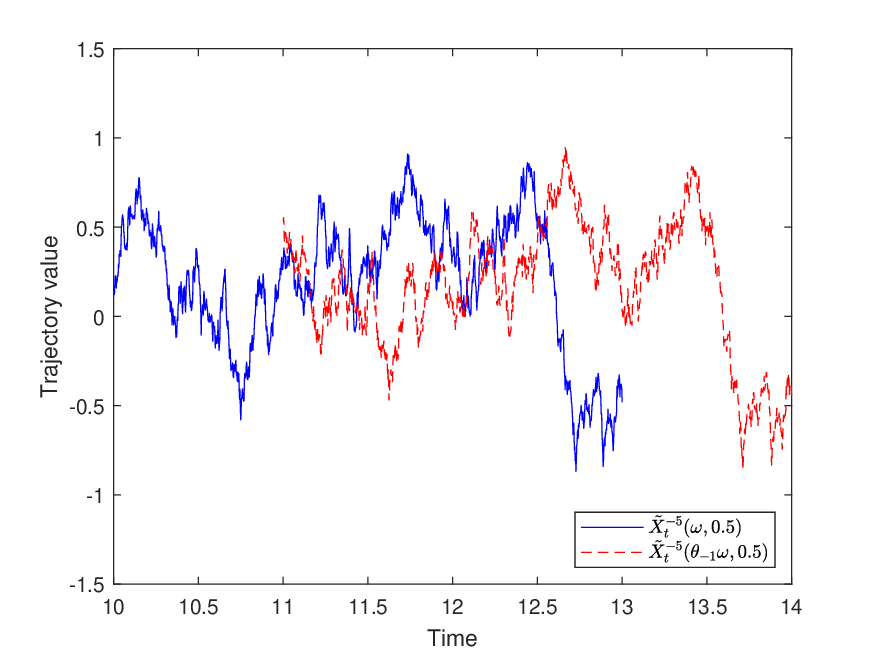}
   \caption[Figure 4.]
	{Simulations of the processes {$\tilde{X}^{-5}_{t}(\omega,0.5),10\leq t\leq 13$} and {$\tilde{X}^{-5}_{t}(\theta_{-1}\omega,0.5),11\leq t\leq 14$}.
	}
 \label{fig4}
\end{figure}

{\color{black}In the second example, we test the performance of the projected Euler method \eqref{eq:the_projected_euler_method} with additive noise as follows:
\begin{equation}
\label{eq_PEM:example_2}
\dd X^{t_{0}}_{t}
=\Big(-\pi X^{t_{0}}_{t}-(X^{t_{0}}_{t})^3
+\sin(2\pi t)\Big)
\,\dd t
+\dd W_{t},
X^{t_0}_{t_0}=\xi.
\end{equation}
One can check that the associated period is 1 and Assumption \ref{ass_PEM} and \ref{ass:add} are fulfilled.}
We conduct a similar experiment to verify the periodicity,
as described in Section \ref{sec:eg1}.
The patterns of $\tilde{X}^{-5}_{t-1}(\omega,0.5)$ 
over $10\leq t \leq 13$ 
and $\tilde{X}^{-5}_{t}(\theta_{-1}\omega,0.5)$ 
over $11 \leq t \leq 14$ 
show identical results under the same realisation $\omega$ in Figure \ref{fig4}.

\begin{figure}[H]
  \centering	\includegraphics[scale=0.5]{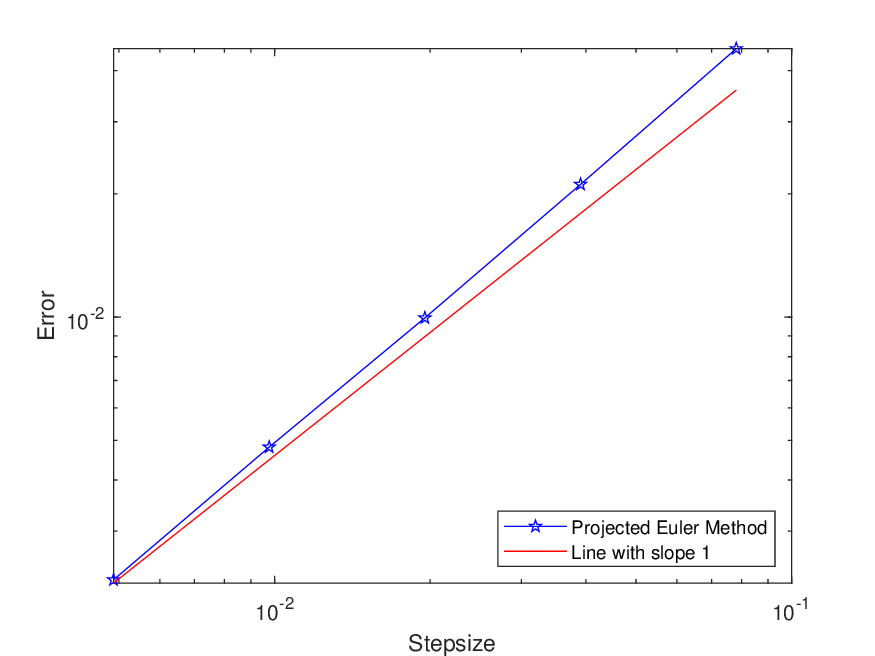}
   \caption[Figure 5.]
	 {
      The mean-square error plot of the projected Euler method \eqref{eq:the_projected_euler_method} for simulating the solution of \eqref{eq_PEM:example_2}.}
 \label{fig5}
\end{figure}

{\color{black}
The performance of the projected Euler method
is also evaluated in terms of mean-square error for simulating SDE \eqref{eq_PEM:example_2} over $[-5, 15]$. 
The comparison between the error line and the reference line in the figure \ref{fig5}
indicates a close match in slopes, supporting an order-one convergence. 
A least squares fit
yields a rate of 1.0751
with a residual of 0.0195 for \eqref{eq_PEM:example_2}. Therefore, the numerical result is in agreement with a strong order of convergence equal to one, as previously indicated in Theorem 
\ref{thm_PEM:add_error analysis}.}

\vskip6mm

\appendix
\section{Proof of Lemma \ref{lem::coupled_monononcity_condition}}

{\color{black}
\begin{proof}
[Proof of Lemma \ref{lem::coupled_monononcity_condition}]
Setting $y=0$ {\color{black} in \eqref{eq:coupled__momotoncity_condition}},
according to \eqref{eq_PEM:coercivity_condition}, we can have,
\begin{equation*}
    \langle 
    x, 
    f(t,x)
    \rangle
    +
    \dfrac{2p_{1}-1}{2}
    \|
    g(t,x)-g(t,0)
    \|^{2}
    \leq 
    \alpha_{1} 
     \| x\|^{2}
     +
     \langle
     x, 
    f(t,0)
    \rangle,
\end{equation*}
then,    
\begin{equation*}
\begin{split}
    &\langle
    x, 
    f(t,x)
    \rangle
    +
    \dfrac{2p_{1}-1}{2}
    \|
    g(t,x)
    \|^{2} \\
    & \quad \leq 
    \alpha_{1} 
     \| x\|^{2}
     +
     \langle
     x, 
    f(t,0)
    \rangle,
    +
    \dfrac{2p_{1}-1}{2}
    \langle
    2g(t,x),
    g(t,0)
    \rangle
    +\dfrac{2p_{1}-1}{2}
    \|g(t,0)\|^2. \\
\end{split} 
\end{equation*}
Using $2ab\leq \epsilon a^2+\tfrac{b^2}{\epsilon}$,
\begin{equation*}
    \langle
     x, 
    f(t,0)
    \rangle
    \leq
    \tfrac{\epsilon}{2}
    \|x\|^2
    +    \tfrac{\|f(t,0)\|^2}{2\epsilon},
\end{equation*}
similarly,
\begin{equation*}
\begin{split}
    \dfrac{2p_{1}-1}{2}
    \langle
    2g(t,x),
    g(t,0)
    \rangle
    &\leq
    \dfrac{2p_{1}-1}{2}
    \times
    \dfrac{2(p_1-p_2)}{2p_{1}-1}
    \|g(t,x)\|^2
    +
    \dfrac{(2p_{1}-1)}{2}
    \times
    \dfrac{2p_{1}-1}{2(p_1-p_2)}\|g(t,0)\|^2 \\
    & =
    \dfrac{2(p_1-p_2)}{2}
    \|g(t,x)\|^2
    +\dfrac{(2p_1-1)^2}{4(p_1-p_2)}
    \|g(t,0)\|^2.
\end{split}
\end{equation*}
Hence,
\begin{equation*}
\begin{split}
    &\langle
    x, 
    f(t,x)
    \rangle
    +
    \dfrac{{\color{black}2p_{2}}-1}{2}
    \|
    g(t,x)
    \|^{2} \\
    & \quad \leq
    (\alpha_1+\epsilon)
    \|x\|^2
    +
    \dfrac{\|f(t,0)\|^2}{2\epsilon}
    +
    \dfrac{(2p_1-1)^2}{4(p_1-p_2)}
    \|g(t,0)\|^2
    +
    \dfrac{2p_{1}-1}{2}
    \|g(t,0)\|^2.
\end{split}    
\end{equation*}
\end{proof}
}
\end{document}